\documentclass[9pt]{article}
\usepackage{mathrsfs}
\usepackage{amsthm}
\usepackage{amssymb}
\usepackage{amsmath}
\usepackage{graphicx}
\usepackage{color}
\usepackage{amsfonts}
\usepackage{float}
\usepackage{cite}
\usepackage [latin1]{inputenc}
\usepackage[text={140mm,210mm},left=35mm,vmarginratio=1:1]{geometry}
\newtheorem{theorem}{Theorem}[section]

\newtheorem{lemma}[theorem]{Lemma}
\newtheorem{proposition}[theorem]{Proposition}

\newtheorem{conjecture}[theorem]{Conjecture}

\numberwithin{equation}{section}
\normalsize

\begin{document}
\title{\textbf{Sample path central limit theorem for the occupation time of the voter model on a lattice}}

\author{Xiaofeng Xue \thanks{\textbf{E-mail}: xfxue@bjtu.edu.cn \textbf{Address}: School of Mathematics and Statistics, Beijing Jiaotong University, Beijing 100044, China.}\\ Beijing Jiaotong University}

\date{}
\maketitle

\noindent {\bf Abstract:} In this paper, we extend the central limit theorem of the occupation time of the voter model on the lattice $\mathbb{Z}^d$ given in \cite{Cox1983} to the sample path case for $d\geq 3$. The proof of our main result utilizes the resolvent strategy and the Poisson flow strategy introduced in previous literatures, where the duality relationship between the voter model and the coalescing random walk plays the key role. For $d=2$ case, we give a conjecture about an analogue result of our main theorem.

\quad



\section{Introduction}\label{section one}

In this paper, we extend the occupation time central limit theorem of the voter model on the lattice $\mathbb{Z}^d$ proved in \cite{Cox1983} to the sample path case when $d\geq 3$. We first introduce some notations for later use. For $x=(x_1,\ldots,x_d)\in \mathbb{Z}^d$, we denote by $\|x\|_1$ the $l_1$-norm of $x$, i.e.,
\[
\|x\|_1=\sum_{j=1}^d|x_j|.
\]
For any $x,y\in \mathbb{Z}^d$, we write $x\sim y$ when and only when $\|x-y\|_1=1$, i.e., $x$ and $y$ are neighbors. We denote by $O$ the origin of $\mathbb{Z}^d$, i.e., $O=(0,0,\ldots,0)$. Now we recall the definition of the voter model. The voter model $\{\eta_t\}_{t\geq 0}$ on $\mathbb{Z}^d$ is a continuous-time Markov process with the state space $\{0, 1\}^{\mathbb{Z}^d}$ and the generator $\mathcal{L}$ given by
\begin{equation}\label{equ 1.1 generator}
\mathcal{L}f(\eta)=\sum_{x\in \mathbb{Z}^d}\sum_{y:y\sim x}\left[f(\eta^{x,y})-f(\eta)\right]
\end{equation}
for any $\eta\in \{0, 1\}^{\mathbb{Z}^d}$ and $f$ from $\{0, 1\}^{\mathbb{Z}^d}$ to $\mathbb{R}$ depending on finite coordinates, where
\[
\eta^{x,y}(z)=
\begin{cases}
\eta(z) & \text{~if~}z\neq x,\\
\eta(y) & \text{~if~} z=x.
\end{cases}
\]
Note that we care about the order of $x$ and $y$, i.e., $\eta^{x,y}\neq \eta^{y,x}$.

Intuitively, the voter model describes the variations of respective supporters of two opposite opinions, which are denoted by $0$ and $1$, of a topic. For each pair of neighbors $x$ and $y$ on $\mathbb{Z}^d$, $x$ adopts $y$'s opinion at rate $1$. For a detailed survey of basic properties of the voter model, see Chapter 5 of \cite{Lig1985} and Part {\rm \uppercase\expandafter{\romannumeral2}} of \cite{Lig1999}.

For any probability measure $\mu$ on $\{0, 1\}^{\mathbb{Z}^d}$, we denote by $\mathbb{P}_\mu$ the probability measure of $\{\eta_t\}_{t\geq 0}$ when initially $\eta_0$ is distributed with $\mu$. We denote by $\mathbb{E}_\mu$ the expectation with respect to $\mathbb{P}_\mu$. For any $\eta\in \{0, 1\}^{\mathbb{Z}^d}$, we write $\mathbb{P}_\mu$ and $\mathbb{E}_\mu$ as $\mathbb{P}_\eta$ and $\mathbb{E}_\eta$ respectively when $\mu$ is the Dirac measure concentrated on $\eta$. For given $0<p<1$, we denote by $\nu_p$ the product measure on $\{0, 1\}^{\mathbb{Z}^d}$ with density $p$, i.e., $\{\eta(x)\}_{x\in \mathbb{Z}^d}$ are independent and identically distributed under $\nu_p$ such that
\[
\nu_p\left(\eta(O)=1\right)=p=1-\nu_p\left(\eta(O)=0\right).
\]
According to the duality relationship between the voter model and the coalescing random walk, which we will recall in Section \ref{section three}, we have
\begin{equation}\label{equ 1.2 density conservation}
\mathbb{E}_{\nu_p}\eta_t(x)=p
\end{equation}
for any $x\in \mathbb{Z}^d$ and $t\geq 0$.

The occupation time process $\{\xi_t\}_{t\geq 0}$ on the origin $O$ is defined as
\[
\xi_t=\int_0^t \eta_s(O)ds.
\]
By Equation \eqref{equ 1.2 density conservation}, when the voter model starts from $\nu_p$, it is reasonable to define the centered occupation time on $O$ as
\[
\int_0^t\left(\eta_s(O)-p\right)ds=\xi_t-tp
\]
for all $t\geq 0$. The following proposition, which is proved in \cite{Cox1983} by Cox and Griffeath, gives the central limit theorem of the aforesaid centered occupation time.
\begin{proposition}\label{proposition occupation time CLT}
(\cite{Cox1983}, Cox and Griffeath, 1983). Let $p\in (0, 1)$ be fixed. If $d\geq 2$ and $\eta_0$ is distributed with $\nu_p$, then
\[
\frac{1}{h_d(t)}\int_0^t\left(\eta_s(O)-p\right)ds
\]
converges weakly to a Gaussian distribution as $t\rightarrow+\infty$, where
\begin{equation}\label{equ hd}
h_d(t)=
\begin{cases}
\frac{t}{\sqrt{\log t}} & \text{~if~}d=2,\\
\\
t^{\frac{3}{4}} & \text{~if~}d=3,\\
\\
\sqrt{t\log t} & \text{~if~}d=4,\\
\\
\sqrt{t} & \text{~if~}d\geq 5.
\end{cases}
\end{equation}
\end{proposition}

In this paper, when $d\geq 3$, we extend the above central limit theorem to the sample path case, i.e., for given $T>0$, we give the weak limit of
\[
\left\{\frac{1}{h_d(N)}\int_0^{tN}\left(\eta_s(O)-p\right)ds:~0\leq t\leq T\right\}_{N\geq 1}
\]
under the Skorohod topology as $N\rightarrow+\infty$. We show that the aforesaid weak limit is a Brownian motion when $d\geq 4$ and a Gaussian process without stationary independent increments when $d=3$. The covariance functions of the weak limit in $d=3$ case are also given. For $d=2$ case, we conjecture that an analogue of our main theorem holds. For the precise statements of our main result and conjecture, see Section \ref{section two}.

Since 1980s, the limit theorem of the occupation time of an interacting particle system is a popular research topic. Central limit theorems, large deviations and moderate deviations are widely discussed for occupation times of models such as exclusion processes (see \cite{Gao2024, Kipnis1987, Landim1992, Lee2004}), voter models (see \cite{Bramson1988, Cox1983, Maillard2009}), contact processes (see \cite{Schonmann1986}), binary contact path processes (see \cite{Xue2024}), branching random walks (see \cite{Birkner2007}), branching Brownian motions (see \cite{Deuschel1998}) and branching $\alpha$-stable processes (see \cite{Li2011}). This paper is inspired a lot by \cite{Kipnis1987} and \cite{Birkner2007}, since the proof of our main result utilizes the resolvent strategy and the Poisson flow strategy introduced in these two literatures. For mathematical details, see Sections \ref{section four d=3} to \ref{section seven d>=5}.

\section{Main result}\label{section two}

In this section, we give our main result. We first introduce some notations for later use. We denote by $\{W_t\}_{t\geq 0}$ the standard Brownian motion with $W_0=0$. We denote by $\{\zeta_t\}_{t\geq 0}$ the Gaussian process with continuous sample path, mean zero and covariance functions given by
\[
{\rm Cov}\left(\zeta_t, \zeta_s\right)=s^{\frac{3}{2}}+t^{\frac{3}{2}}-\frac{1}{2}(t-s)^{\frac{3}{2}}-\frac{1}{2}(t+s)^{\frac{3}{2}}
\]
for any $0\leq s\leq t$.


Throughout this paper, we let $p\in (0, 1)$ be a fixed constant and $T>0$ be a fixed moment. We let $\nu_p$ be the product measure on $\{0, 1\}^{\mathbb{Z}^d}$ with density $p$ defined as in Section \ref{section one} and denote by $C[0, T]$ the set of continuous functions from $[0, T]$ to $\mathbb{R}$. We let $C[0, T]$ be endowed with the Skorohod topology. We denote by $\{X_t\}_{t\geq 0}$ the continuous-time simple random walk on $\mathbb{Z}^d$ with generator $\mathcal{G}$ given by
\[
\mathcal{G}F(x)=\sum_{y:y\sim x}\left[F(y)-F(x)\right]
\]
for any $x\in \mathbb{Z}^d$ and bounded $F$ from $\mathbb{Z}^d$ to $\mathbb{R}$. We denote by $p_t(\cdot, \cdot)$ the transition probabilities of $\{X_t\}_{t\geq 0}$, i.e, for any $x,y\in \mathbb{Z}^d$,
\[
p_t(x, y)=\mathbb{P}\left(X_t=y\big|X_0=x\right).
\]
For $d\geq 3$, we use $\gamma$ to denote
\[
\mathbb{P}\left(X_t\neq O\text{~for all~}t\geq 0\Bigg|\|X_0\|_1=1\right).
\]
According to the strong Markov property of $\{X_t\}_{t\geq 0}$ and the fact that $\{X_t\}_{t\geq 0}$ stays at each state for an exponential time with rate $2d$, we have
\begin{equation}\label{equ return green}
\int_0^{+\infty}p_s(O,O)ds=\frac{1}{2d}\frac{1}{\gamma}.
\end{equation}
We write $p_t(\cdot, \cdot), \gamma$ as $p_{t,d}(\cdot, \cdot), \gamma_d$ when we need to emphasize the dependence on $d$. For $d\geq 3$, we define $C_d$ as
\[
C_d=
\begin{cases}
\sqrt{\frac{4p(1-p)\gamma_3}{\pi^{\frac{3}{2}}}} & \text{~if~}d=3,\\
\\
\sqrt{\frac{\gamma_4p(1-p)}{\pi^2}} & \text{~if~}d=4,\\
\\
\sqrt{4dp(1-p)\gamma_d\int_0^{+\infty}sp_{s,d}(O, O)ds} & \text{~if~}d\geq 5.
\end{cases}
\]
Note that, according to the well-known result that $p_t(O,O)=\Theta(t^{-\frac{d}{2}})$ (see Chapter 2 of \cite{Lawler2010}),
\[
\int_0^{+\infty}sp_{s, d}(O,O)ds<+\infty
\]
when and only when $d\geq 5$.

Now we give our main result.
\begin{theorem}\label{theorem 2.1 main sample path CLT}
Let $d\geq 3$ and $h_d(\cdot)$ be defined as in \eqref{equ hd}. If $\eta_0$ is distributed with $\nu_p$, then $\left\{\frac{1}{h_d(N)}\int_0^{tN}\left(\eta_s(O)-p\right)ds:~0\leq t\leq T\right\}_{N\geq 1}$ converges weakly, with respect to the Skorohod topology of $C[0, T]$, to $\{\Lambda_t\}_{0\leq t\leq T}$ as $N\rightarrow+\infty$, where
\[
\{\Lambda_t\}_{0\leq t\leq T}=
\begin{cases}
\left\{C_3\zeta_t\right\}_{0\leq t\leq T} & \text{~if~} d=3,\\
\\
\left\{C_dW_t\right\}_{0\leq t\leq T} & \text{~if~}d\geq 4.
\end{cases}
\]

\end{theorem}

The proof of Theorem \ref{theorem 2.1 main sample path CLT} is inspired by the resolvent strategy introduced in \cite{Kipnis1987} and the Poisson flow strategy introduced in \cite{Birkner2007}. By utilizing the resolvent strategy, we can decompose the centered occupation time as a martingale plus an error term which converges weakly to $0$. Inspired by the Poisson flow strategy, we write $\{\eta_t(x):~x\in \mathbb{Z}^d, t\geq 0\}$ as the solutions to a series of stochastic differential equations driven by a family of independent Poisson processes and then we can show that the aforesaid martingale part converges weakly to the target Gaussian process by applying the properties of Poisson processes. For the mathematical details, see Sections \ref{section four d=3}-\ref{section seven d>=5}.

We denote by $\{\vartheta_t\}_{t\geq 0}$ the Gaussian process with continuous sample path, mean zero and covariance functions given by
\[
{\rm Cov}\left(\vartheta_t, \vartheta_s\right)=\frac{(t+s)^2}{4}\log(t+s)+\frac{(t-s)^2}{4}\log(t-s)-\frac{s^2\log s}{2}-\frac{t^2\log t}{2}
\]
for any $0\leq s\leq t$. About the $d=2$ case, we have the following conjecture.

\begin{conjecture}\label{conjecture in d=2 case}
Let $d=2$ and $\eta_0$ be distributed with $\nu_p$,  then
\[
\left\{\left(\frac{N}{\sqrt{\log N}}\right)^{-1}\int_0^{tN}\left(\eta_s(O)-p\right)ds:~0\leq t\leq T\right\}_{N\geq 1}
\]
converges weakly, with respect to the Skorohod topology of $C[0, T]$, to $\{C_2\vartheta_t\}_{0\leq t\leq T}$ as $N\rightarrow+\infty$, where $C_2=\sqrt{2p(1-p)}$.

\end{conjecture}

In Section \ref{section d=2 conjecture comment}, we explain why we believe Conjecture \ref{conjecture in d=2 case} holds. We will show that what we lack to prove Conjecture \ref{conjecture in d=2 case} rigorously is an analogue of a martingale variance control given in the proof of Theorem \ref{theorem 2.1 main sample path CLT} in the $d=3$ case. For mathematical details, see Section \ref{section d=2 conjecture comment}. 

\section{Duality relationships}\label{section three}
\begin{figure}[H]
\centering
\includegraphics[scale=0.3]{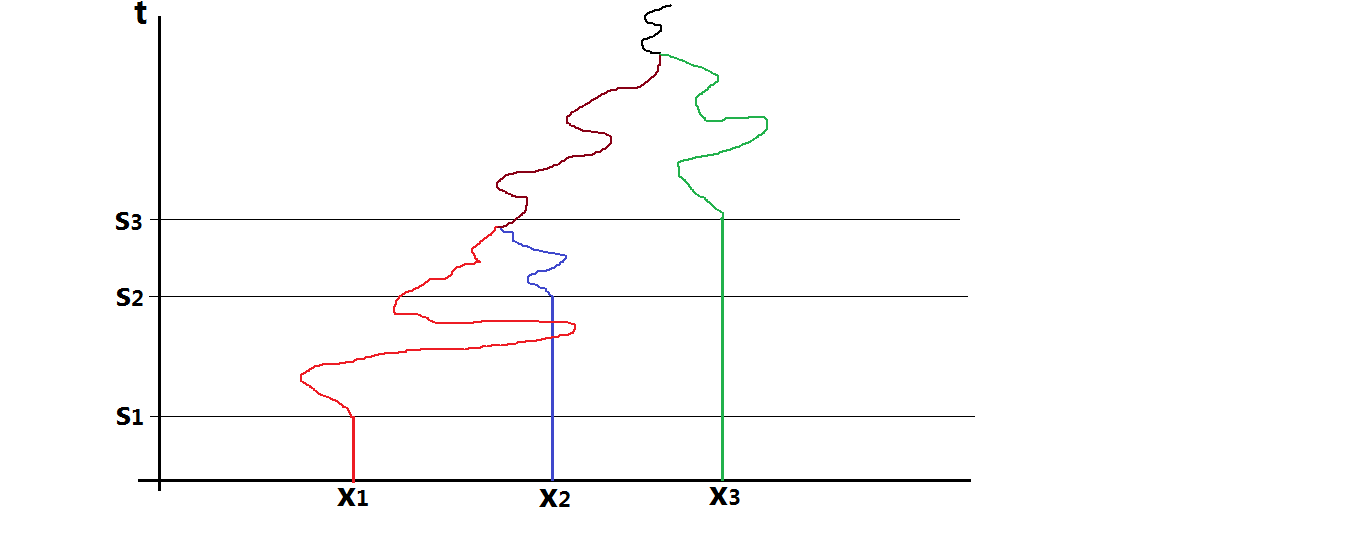}
\caption{A coalescing random walk example}\label{figure coalescing random walk}
\end{figure}

As a preliminary for the proof of Theorem \ref{theorem 2.1 main sample path CLT}, in this section we recall the duality relationship between the voter model and the coalescing random walks. For given integer $m\geq 1$ and $0\leq s_1\leq s_2\leq \ldots\leq s_m$, $x_1, x_2,\ldots,x_m\in \mathbb{Z}^d$, the coalescing random walk
\[
\left\{\left(\mathcal{X}_{t,1}^{x_1, s_1}, \mathcal{X}_{t,2}^{x_2,s_2},\ldots, \mathcal{X}_{t, m}^{x_m, s_m}\right)\right\}_{t\geq 0}
\]
evolves as follows. For all $1\leq k\leq m$, $\mathcal{X}_{t,k}^{x_k, s_k}=x_k$ for $0\leq t\leq s_k$ and $\{\mathcal{X}_{s_k+t,k}^{x_k, s_k}\}_{t\geq 0}$ is a copy of $\{X_t\}_{t\geq 0}$ starting from $x_k$, where $X_t$ is defined as in Section \ref{section two}. We call $\mathcal{X}_{t, k}^{x_k, s_k}$ active when $t\geq s_k$. All active walks evolve independently until collisions occur. When two independent active walks collide with each other, then they are coalesced as one active walk. Figure \ref{figure coalescing random walk} gives a picture describing an example of the coalescing random walk with $m=3$.

We denote by $\hat{\mathbb{P}}$ the probability measure of the coalescing random walk and by $\hat{\mathbb{E}}$ the expectation with respect to $\hat{\mathbb{P}}$. The following proposition gives the duality relationship between the voter model and the coalescing random walk,

\begin{proposition}\label{proposition duality}
For any integer $m\geq 1$, $y_1, y_2,\ldots, y_m\in \mathbb{Z}^d$, $0\leq t_1\leq t_2\leq \ldots\leq t_m$, $\eta\in \{0, 1\}^{\mathbb{Z}^d}$ and function $G$ from $\{0, 1\}^m$ to $\mathbb{R}$, we have
\begin{align}\label{equ duality relationship}
&\mathbb{E}_\eta G\left(\eta_{t_1}(y_1), \eta_{t_2}(y_2),\ldots, \eta_{t_m}(y_m)\right)\\
&=\hat{\mathbb{E}}G\left(\eta\left(\mathcal{X}_{t_m, m}^{x_m, s_m}\right), \eta\left(\mathcal{X}_{t_m, m-1}^{x_{m-1}, s_{m-1}}\right),\ldots, \eta\left(\mathcal{X}_{t_m, 2}^{x_2, s_2}\right), \eta\left(\mathcal{X}_{t_m, 1}^{x_1, s_1}\right)\right),\notag
\end{align}
where
\[
\left(x_1,x_2,\ldots,x_m\right)=\left(y_m, y_{m-1},\ldots,y_1\right)
\]
and
\[
\left(s_1,s_2,\ldots,s_m\right)=\left(t_m-t_m, t_m-t_{m-1}, \ldots, t_m-t_2, t_m-t_1\right).
\]
Consequently, for given $0<p<1$,
\begin{align}\label{equ duality product}
&\mathbb{E}_{\nu_p} G\left(\eta_{t_1}(y_1), \eta_{t_2}(y_2),\ldots, \eta_{t_m}(y_m)\right)\\
&=\int_{\{0, 1\}^{\mathbb{Z}^d}}\left(\hat{\mathbb{E}}G\left(\eta\left(\mathcal{X}_{t_m, m}^{x_m, s_m}\right), \eta\left(\mathcal{X}_{t_m, m-1}^{x_{m-1}, s_{m-1}}\right),\ldots, \eta\left(\mathcal{X}_{t_m, 2}^{x_2, s_2}\right), \eta\left(\mathcal{X}_{t_m, 1}^{x_1, s_1}\right)\right)\right)\nu_p(d\eta).\notag
\end{align}
\end{proposition}

The proof of Proposition \ref{proposition duality} can be found in Chapter 3 of \cite{Lig1985}, which utilizes a graphical method introduced in \cite{Har1978}.

For later use, here we give some applications of Proposition \ref{proposition duality}. Let $m=1$ and $G(u)=u$ for $u=0, 1$, then for $x\in \mathbb{Z}^d$ we have
\begin{equation}\label{equ dual eta(x)}
\mathbb{P}_\eta\left(\eta_t(x)=1\right)=\hat{\mathbb{P}}\left(\eta(\mathcal{X}_{t,1}^{x,0})=1\right)=\sum_{y\in \mathbb{Z}^d}p_t(x,y)\eta(y)
\end{equation}
and hence
\begin{equation}\label{equ conservative mean}
\mathbb{E}_{\nu_p}\eta_t(x)=\mathbb{P}_{\nu_p}(\eta_t(x)=1)=\sum_{y\in \mathbb{Z}^d}p_t(x,y)p=p
\end{equation}
as we have introduced in Section \ref{section one}.

Let $m=2$ and $G(u, v)=(u-v)^2=1_{\{u\neq v\}}$ for $u,v\in \{0, 1\}$, then for different $x,y\in \mathbb{Z}^d$ we have
\begin{align}\label{equ duality minus square}
\mathbb{E}_{\nu_p}\left((\eta_t(x)-\eta_t(y))^2\right)&=\int_{\{0,1\}^{\mathbb{Z}^d}}\hat{\mathbb{P}}\left(\eta(\mathcal{X}_{t,1}^{x,0})\neq \eta(\mathcal{X}_{t,2}^{y, 0})\right)\nu_p(d\eta) \notag\\
&=\int_{\{0,1\}^{\mathbb{Z}^d}}\hat{\mathbb{P}}\left(\eta(\mathcal{X}_{t,1}^{x,0})\neq \eta(\mathcal{X}_{t,2}^{y, 0}), \tau_{x,y}>t\right)\nu_p(d\eta) \\
&=2p(1-p)\hat{\mathbb{P}}(\tau_{xy}>t), \notag
\end{align}
where
\[
\tau_{xy}=\inf\left\{t:~\mathcal{X}_{t,1}^{x,0}=\mathcal{X}_{t,2}^{y, 0}\right\}.
\]
Note that the last equation in \eqref{equ duality minus square} utilizes the fact that
\[
\mathbb{P}_{\nu_p}(\eta(z)\neq \eta(w))=\mathbb{P}_{\nu_p}\left(\eta(z)=1, \eta(w)=0\right)+\mathbb{P}_{\nu_p}\left(\eta(z)=0, \eta(w)=1\right)=2p(1-p)
\]
for any different $z,w\in \{0, 1\}^{\mathbb{Z}^d}$.

Let $m=2$ and $G(u, v)=(u-p)(v-p)$ for $u, v\in \{0, 1\}$, then for $x, y\in \mathbb{Z}^d$ and $t\geq 0$ we have
\begin{align}\label{equ voter covariance}
{\rm Cov}_{\nu_p}(\eta_t(x), \eta_t(y))&=\mathbb{E}_{\nu_p}G(\eta_t(x), \eta_t(y)) \notag\\
&=\hat{\mathbb{E}}\left(\int_{\{0, 1\}^{\mathbb{Z}^d}}\left(\eta(\mathcal{X}_{t,1}^{x,0})-p\right)\left(\eta(\mathcal{X}_{t,2}^{y, 0})-p\right)\nu_p(d \eta)\right) \notag\\
&=\left((1-p)^2p+(1-p)(-p)^2\right)\hat{\mathbb{P}}\left(\tau_{xy}\leq t\right)\notag\\
&=(1-p)p\hat{\mathbb{P}}\left(\tau_{xy}\leq t\right).
\end{align}
Note that the third equation in the above formula utilizes the fact that
\[
\int_{\{0, 1\}^{\mathbb{Z}^d}}\left(\eta(v)-p\right)\left(\eta(u)-p\right)\nu_p(d \eta)=0
\] for $u\neq v$.

Let $m=4$ and $G(u,v,a,b)=(u-p)(v-p)(a-p)(b-p)$ for $u,v,a,b\in \{0, 1\}$, then for four different $x_1, x_2, x_3, x_4\in \mathbb{Z}^d$ and $0\leq t_1\leq t_2\leq t_3\leq t_4$, we have
\begin{align}\label{equ 3.4}
\mathbb{E}_{\nu_p}\left(\prod_{j=1}^4\left(\eta_{t_j}(x_j)-p\right)\right)=\int_{\{0, 1\}^{\mathbb{Z}^d}}\hat{\mathbb{E}}\left(\prod_{j=1}^4\left(\eta(\mathcal{X}_{t_4, j}^{x_j,t_4-t_j})-p\right)\right)\nu_p(d\eta).
\end{align}
For $1\leq i\neq j\leq 4$, we define
\[
\tau_{x_ix_j}^{t_4-t_i, t_4-t_j}=\inf\left\{t:~t>t_4-t_i, t>t_4-t_j, \mathcal{X}_{t, i}^{x_i,t_4-t_i}=\mathcal{X}_{t, j}^{x_j,t_4-t_j}\right\}.
\]
For given $z_1,z_2,z_3,z_4\in \mathbb{Z}^d$, it is easy to check that
\begin{align}\label{equ 3.5}
&\int_{\{0, 1\}^{\mathbb{Z}^d}}\prod_{j=1}^4\left(\eta(z_j)-p\right)\nu_p(d\eta)\\
&=
\begin{cases}
p(1-p)^4+(1-p)(-p)^4 & \text{~if~}z_1=z_2=z_3=z_4,\\
\left(p(1-p)^2+(1-p)(-p)^2\right)^2 & \text{~if~}z_1=z_2\neq z_3=z_4,\\
\left(p(1-p)^2+(1-p)(-p)^2\right)^2 & \text{~if~}z_1=z_3\neq z_2=z_4,\\
\left(p(1-p)^2+(1-p)(-p)^2\right)^2 & \text{~if~}z_1=z_4\neq z_2=z_3,\\
0 & \text{~else.}
\end{cases}
\notag
\end{align}
For any $1\leq i\neq j\leq 4$, we denote by $\mathcal{A}_{ij}$ the event that
\[
\tau_{x_lx_k}^{t_4-t_l, t_4-t_k}\leq t_4
\]
for both $\{l,k\}=\{i,j\}$ and $\{l,k\}=\{1,2,3,4\}\setminus \{i,j\}$. By Equations \eqref{equ 3.4} and \eqref{equ 3.5}, we have
\begin{equation}\label{equ 3.6}
\mathbb{E}_{\nu_p}\left(\prod_{j=1}^4\left(\eta_{t_j}(x_j)-p\right)\right)\leq \hat{\mathbb{P}}\left(\mathcal{A}_{12}\bigcup\mathcal{A}_{13}\bigcup\mathcal{A}_{14}\right).
\end{equation}
To further bound $\hat{\mathbb{P}}\left(\mathcal{A}_{12}\bigcup\mathcal{A}_{13}\bigcup\mathcal{A}_{14}\right)$ from above, we define $\{Z_{t,j}:~t\geq 0\}_{j=1}^4$ as four independent random walks on $\mathbb{Z}^d$ such that, for all $1\leq j\leq 4$, $Z_{t,j}=x_j$ for all $0\leq t\leq t_4-t_j$ and $\{Z_{t_4-t_j+t,j}\}_{t\geq 0}$ is a copy of $\{X_t\}_{t\geq 0}$ starting from $x_j$. For any $1\leq i\neq j\leq 4$, we define
\[
\hat{\tau}_{ij}=\inf\left\{t:~t>t_4-t_i, t>t_4-t_j, Z_{t, i}=Z_{t, j}\right\}.
\]
Hence, $\hat{\tau}_{ij}$ and $\hat{\tau}_{lk}$ are independent when $\{l,k\}=\{1,2,3,4\}\setminus \{i,j\}$. Now we equivalently define the coalescing random walk as follows. If
\[
\inf\{\hat{\tau}_{ij}:~1\leq i\neq j\leq 4\}=+\infty,
\]
then $\mathcal{X}_{t, j}^{x_j,t_4-t_j}=Z_{t,j}$ for all $1\leq j\leq 4$ and $t\geq 0$. Otherwise, if
\[
\hat{\tau}_{ij}=\inf\{\hat{\tau}_{kl}:~1\leq k\neq l\leq 4\}<+\infty
\]
for some $1\leq i<j\leq 4$, then for $l,k$ such that $l<k$ and $\{l,k\}=\{1,2,3,4\}\setminus\{i, j\}$, we define
$\mathcal{X}_{t, k}^{x_k,t_4-t_k}=Z_{t,k}$ for all $t\geq 0$ and
\[
\mathcal{X}_{t,l}^{x_l,t_4-t_l}=
\begin{cases}
Z_{t,l} \text{~if~}t\leq \hat{\tau}_{lk},\\
Z_{t,k} \text{~if~}t>\hat{\tau}_{lk}.
\end{cases}
\]
Let $\tilde{\tau}_1=\inf\{s:~s>t_4-t_j, s>t_4-t_l, Z_{s,j}=\mathcal{X}_{s,l}^{x_l,t_4-t_l}\}$, $\tilde{\tau}_2=\inf\{s:~s>t_4-t_j, s>t_4-t_k, Z_{s,j}=\mathcal{X}_{s,k}^{x_k,t_4-t_k}\}$ and $\tilde{\tau}=\min\{\tilde{\tau}_1, \tilde{\tau}_2\}$, then
we further define
\[
\mathcal{X}_{t,j}^{x_j, t_4-t_j}=
\begin{cases}
Z_{t,j} &\text{~if~}t\leq \tilde{\tau},\\
\mathcal{X}_{t, l}^{x_l,t_4-t_l} &\text{~if~}t>\tilde{\tau} \text{~and~}\tilde{\tau}=\tilde{\tau}_1,\\
\mathcal{X}_{t, k}^{x_k,t_4-t_k} &\text{~if~}t>\tilde{\tau} \text{~and~}\tilde{\tau}=\tilde{\tau}_2.
\end{cases}
\]
At last, we define
\[
\mathcal{X}_{t,i}^{x_i, t_4-t_i}=
\begin{cases}
Z_{t, i} & \text{~if~}t\leq \hat{\tau}_{ij},\\
\mathcal{X}_{t,j}^{x_j, t_4-t_j} & \text{~if~}t>\hat{\tau}_{ij}.
\end{cases}
\]
It is not difficult to check that the above definition of
\[
\left\{\left(\mathcal{X}_{t, 1}^{x_1,t_4-t_1}, \mathcal{X}_{t, 2}^{x_2,t_4-t_2}, \mathcal{X}_{t, 3}^{x_3,t_4-t_3}, \mathcal{X}_{t, 4}^{x_1,t_4-t_4}\right)\right\}_{t\geq 0}
\]
is equivalent to the one given before Proposition \ref{proposition duality}.

For $1\leq i\neq j\leq 4$, we denote by $\mathcal{B}_{ij}$ the event that
\[
\hat{\tau}_{lk}\leq t_4
\]
for both $\{l,k\}=\{i,j\}$ and $\{l,k\}=\{1,2,3,4\}\setminus \{i,j\}$. According to our second definition of the coalescing random walk, on the event
$\mathcal{A}_{12}\bigcup\mathcal{A}_{13}\bigcup\mathcal{A}_{14}$, there exist $1\leq i<j\leq 4$ such that
\[
\tau_{x_ix_j}^{t_4-t_i, t_4-t_j}=\hat{\tau}_{ij}=\inf\{\hat{\tau}_{kl}:~1\leq k\neq l\leq 4\}<t_4
\]
and
\[
\hat{\tau}_{lk}=\tau_{x_lx_k}^{t_4-t_l, t_4-t_k}<t_4
\]
for $\{l,k\}=\{1,2,3,4\}\setminus\{i, j\}$. As a result, the event $\mathcal{B}_{12}\bigcup\mathcal{B}_{13}\bigcup\mathcal{B}_{14}$ occurs. Consequently,
\begin{align}\label{equ 3.7}
&\hat{\mathbb{P}}\left(\mathcal{A}_{12}\bigcup\mathcal{A}_{13}\bigcup\mathcal{A}_{14}\right)\leq \hat{\mathbb{P}}\left(\mathcal{B}_{12}\bigcup\mathcal{B}_{13}\bigcup\mathcal{B}_{14}\right) \notag\\
&\leq \hat{\mathbb{P}}\left(\hat{\tau}_{12}\leq t_4\right)\hat{\mathbb{P}}\left(\hat{\tau}_{34}\leq t_4\right)
+\hat{\mathbb{P}}\left(\hat{\tau}_{13}\leq t_4\right)\hat{\mathbb{P}}\left(\hat{\tau}_{24}\leq t_4\right)
+\hat{\mathbb{P}}\left(\hat{\tau}_{14}\leq t_4\right)\hat{\mathbb{P}}\left(\hat{\tau}_{23}\leq t_4\right).
\end{align}
Not that the last inequality in \eqref{equ 3.7} utilizes the fact that $\hat{\tau}_{ij}$ and $\hat{\tau}_{lk}$ are independent when $\{l,k\}=\{1,2,3,4\}\setminus \{i,j\}$.

Equation \eqref{equ 3.7} plays the key role in the check of the tightness of the sample path of the centered occupation time in cases $d\geq 4$. For mathematical details, see Sections \ref{section six d=4} and \ref{section seven d>=5}.

\section{The proof of Theorem \ref{theorem 2.1 main sample path CLT}: $d=3$ case}\label{section four d=3}
In this section, we prove Theorem \ref{theorem 2.1 main sample path CLT} in the case $d=3$. As we have introduced in Section \ref{section two}, our proof is inspired by the resolvent strategy introduced in \cite{Kipnis1987} and the Poisson flow strategy introduced in \cite{Birkner2007}. So we first give an equivalent definition of the voter model $\{\eta_t\}_{t\geq 0}$ via a family of independent Poisson processes. For any $x,y\in \mathbb{Z}^3$ such that $x\sim y$, let $\{N_t^{x,y}\}_{t\geq 0}$ be a Poisson process with rate $1$. We assume that all these Poisson processes are independent. Note that we care about the order of $x$ and $y$, i.e., $N_t^{x,y}\neq N_t^{y,x}$. Let $\{\eta_0(x)\}_{x\in \mathbb{Z}^d}$ be distributed with some initial distribution on $\{0,1\}^{\mathbb{Z}^3}$, then, for any $x\in \mathbb{Z}^d$, the value of $\eta_\cdot(x)$ flips only at event moments of $\{N_t^{x,y}:~t\geq 0\}_{y:y\sim x}$. In detail, if $t$ is an event moment of $N^{x,y}_\cdot$, then $\eta_t(x)$ takes $\eta_{t-}(y)$, where $\eta_{t-}(y)=\lim_{s<t,s\uparrow t}\eta_s(y)$. Hence, at the event moment $t$ of the Poisson process $N_\cdot^{x,y}$, the increment of $\eta_\cdot(x)$ is $\eta_{t-}(y)-\eta_{t-}(x)$. As a result, we can consider $\{\eta_t(x):~t\geq 0\}_{x\in \mathbb{Z}^3}$ as the solutions to the following series of stochastic differential equations driven by $\{N^{x,y}_\cdot\}_{x\in \mathbb{Z}^3, y\sim x}$: for all $x\in \mathbb{Z}^3$ and $t\geq 0$,
\begin{equation}\label{equ Poisson flow representation of voter model}
\eta_t(x)=\eta_0(x)+\sum_{y:y\sim x}\int_0^t\left(\eta_{s-}(y)-\eta_{s-}(x)\right)dN_s^{x,y}.
\end{equation}
We further define $\hat{N}_t^{x,y}=N_t^{x,y}-t$, then $\{\hat{N}_t^{x,y}\}_{t\geq 0}$ is a martingale with mean zero and cross variation processes given by
\begin{equation}\label{equ 4.cross variation}
\langle \hat{N}^{x,y}, \hat{N}^{z,w}\rangle_t=
\begin{cases}
t & \text{~if~} (x,y)=(z,w)\\
0 & \text{~else}.
\end{cases}
\end{equation}
For any $f$ from $[0, +\infty)\times\mathbb{Z}^3$ to $\mathbb{R}$ such that $\sum_{x\in \mathbb{Z}^d}|f(t, x)|<+\infty$ for any $t\geq 0$ and $f(\cdot,x)\in C^1[0, +\infty)$ for all $x\in \mathbb{Z}^3$, we define
\[
H_f(t, \eta)=\sum_{x\in \mathbb{Z}^d}(\eta(x)-p)f(t, x)
\]
for any $\eta\in \{0, 1\}^{\mathbb{Z}^3}$ and $t\geq 0$. According to the Dynkin's martingale formula,
\[
\left\{H_f(t, \eta_t)-H_f(0, \eta_0)-\int_0^t\left(\mathcal{L}+\partial_s\right)H_f(s, \eta_s)ds\right\}_{t\geq 0}
\]
is a martingale, where $\mathcal{L}$ is the generator given in \eqref{equ 1.1 generator}. We denote by $M_f(t)$ the martingale $H_f(t, \eta_t)-H_f(0, \eta_0)-\int_0^t\left(\mathcal{L}+\partial_s\right)H_f(s, \eta_s)ds$. We have the following representation lemma.

\begin{lemma}\label{lemma martingale representation}
For any $t\geq 0$,
\[
M_f(t)=\sum_{x\in \mathbb{Z}^3}\sum_{y:y\sim x}\int_0^tf(s, x)\left(\eta_{s-}(y)-\eta_{s-}(x)\right)d\hat{N}_s^{x,y}.
\]
\end{lemma}

According to \eqref{equ Poisson flow representation of voter model} and the fact that
\[
\mathcal{L}H_f(t, \eta)=\sum_{x\in \mathbb{Z}^3}\sum_{y:y\sim x}(\eta(y)-\eta(x))f(t,x),
\]
it is not difficult to prove Lemma \ref{lemma martingale representation} via It\^{o}'s formula and hence we omit details here.

To utilize the resolvent strategy, we define
\[
v(t,x)=\int_0^tp_s(O, x)ds
\]
for any $x\in \mathbb{Z}^3, t\geq 0$ and define
\[
V_s^t(\eta)=\sum_{x\in \mathbb{Z}^3}(\eta(x)-p)v(t-s, x)
\]
for any $t\geq 0, 0\leq s\leq t$ and $\eta\in \{0, 1\}^{\mathbb{Z}^3}$. For $0\leq s\leq t$, we further define
\[
\mathcal{M}_s^t=V_s^t(\eta_s)-V_0^t(\eta_0)-\int_0^s(\partial_r+\mathcal{L})V_r^t(\eta_r)dr,
\]
then $\{\mathcal{M}_s^t\}_{0\leq s\leq t}$ is a martingale and
\begin{equation}\label{equ 4.2}
\mathcal{M}_s^t=\sum_{x\in \mathbb{Z}^3}\sum_{y:y\sim x}\int_0^sv(t-r, x)\left(\eta_{r-}(y)-\eta_{r-}(x)\right)d\hat{N}_r^{x,y}
\end{equation}
by Lemma \ref{lemma martingale representation}. According to the expressions of $\mathcal{L}$ and $v(\cdot, \cdot)$, we have
\begin{align}\label{equ generator and occupation time}
\mathcal{L}V_r^t(\eta_r)&=\sum_{x\in \mathbb{Z}^3}\sum_{y:y\sim x}(\eta_r(y)-\eta_r(x))v(t-r,x)\notag\\
&=\sum_{x\in \mathbb{Z}^3}\eta_r(x)\left(\sum_{y:y\sim x}v(t-r, y)-6v(t-r, x)\right)=\sum_{x\in \mathbb{Z}^3}\eta_r(x)\int_0^{t-r}\frac{d}{d\theta}p_\theta(O, x)d\theta\notag\\
&=\sum_{x\in \mathbb{Z}^3}\eta_r(x)\left(p_{t-r}(O, x)-1_{\{x=O\}}\right)=\sum_{x\in Z^3}(\eta_r(x)-p)p_{t-r}(O, x)-(\eta_r(O)-p)\notag\\
&=-\partial_rV_r^t(\eta_r)-(\eta_r(O)-p).
\end{align}
Therefore,
\[
\mathcal{M}_s^t=V_s^t(\eta_s)-V_0^t(\eta_0)+\int_0^s\left(\eta_r(O)-p\right)dr
\]
and hence
\begin{equation}\label{equ 4.3}
\frac{1}{N^{\frac{3}{4}}}\int_0^{tN}(\eta_s(O)-p)ds=\frac{1}{N^{\frac{3}{4}}}\mathcal{M}_{tN}^{tN}+\frac{1}{N^{\frac{3}{4}}}V_0^{tN}(\eta_0)
\end{equation}
since $V_t^t(\eta)=0$. The following lemma shows that $\frac{1}{N^{\frac{3}{4}}}V_0^{tN}(\eta_0)$ is a small error term.
\begin{lemma}\label{lemma small error}
For any given $t>0$, if $\eta_0$ is distributed with $\nu_p$, then
\[
\lim_{N\rightarrow+\infty}\frac{1}{N^{\frac{3}{4}}}V_0^{tN}(\eta_0)=0
\]
in $L^2$.
\end{lemma}

\proof
According to the definition of $V_0^t$, we have
\[
\mathbb{E}_{\nu_p}\left(\left(\frac{1}{N^{\frac{3}{4}}}V_0^{tN}(\eta_0)\right)^2\right)=\frac{p(1-p)}{N^{\frac{3}{2}}}\int_0^{tN}\int_0^{tN}p_{\theta+r}(O, O)d\theta dr.
\]
According to the local central limit theorem of the simple random walk on $\mathbb{Z}^d$ (see Chapter 2 of \cite{Lawler2010}), we have $p_t(O, O)=O(t^{-\frac{3}{2}})$ and consequently
\[
\mathbb{E}_{\nu_p}\left(\left(\frac{1}{N^{\frac{3}{4}}}V_0^{tN}(\eta_0)\right)^2\right)=\frac{O(N^{\frac{1}{2}})}{N^{\frac{3}{2}}}=O(N^{-1}),
\]
which completes this proof.

\qed

The following lemma shows that any finite dimensional distribution of $\{\frac{1}{N^{\frac{3}{4}}}\mathcal{M}_{tN}^{tN}\}_{0\leq t\leq T}$ converges weakly to a multivariate Gaussian distribution.

\begin{lemma}\label{lemma 4.3}
Let $\eta_0$ be distributed with $\nu_p$. For any integer $m\geq 1$ and $0\leq t_1<t_2<\ldots<t_m\leq T$, $\left(\frac{1}{N^{\frac{3}{4}}}\mathcal{M}_{t_1N}^{t_1N}, \frac{1}{N^{\frac{3}{4}}}\mathcal{M}_{t_2N}^{t_2N}, \ldots, \frac{1}{N^{\frac{3}{4}}}\mathcal{M}_{t_mN}^{t_mN}\right)$ converges weakly, as $N\rightarrow+\infty$, to a $\mathbb{R}^m$-valued Gaussian random variable $(\Xi_1,\ldots,\Xi_m)$, which depends on the choice of $t_1, t_2,\ldots, t_m$. Furthermore,
\[
\lim_{N\rightarrow+\infty}\mathbb{E}_{\nu_p}\left(\left(\sum_{j=1}^mK_j\frac{1}{N^{\frac{3}{4}}}\mathcal{M}_{t_jN}^{t_jN}\right)^2\right)
=\mathbb{E}\left(\left(\sum_{j=1}^mK_j\Xi_j\right)^2\right)
\]
for any constants $K_1,\ldots,K_m$.
\end{lemma}

As a preliminary of the proof of Lemma \ref{lemma 4.3}, we introduce some notations and definitions. For any $u\in \mathbb{R}^3$ and integer $N\geq 1$, we denote by $u_N$ the vector in $\mathbb{R}^3$ such that $u_N\in \frac{\mathbb{Z}^3}{\sqrt{N}}$ and $u-u_N\in (-\frac{1}{2\sqrt{N}}, \frac{1}{2\sqrt{N}}]^3$.

For each $N\geq 1$, we denote by $\mathcal{Y}^N$ the random measure on $[0, T]\times \mathbb{R}^3$ such that
\begin{align*}
\mathcal{Y}^N(H)&=N^{\frac{1}{4}}\int_{\mathbb{R}^3}\sum_{y:y\sim \sqrt{N}u_N}\left(\int_0^TH(s,u)\left(\eta_{Ns-}(y)-\eta_{Ns-}(\sqrt{N}u_N)\right)d\hat{N}_{Ns}^{\sqrt{N}u_N,y}\right)du\\
&=N^{\frac{1}{4}}\sum_{x\in \mathbb{Z}^3}\sum_{y:y\sim x}\int_0^T\left(\int_{\frac{x}{\sqrt{N}}+(-\frac{1}{2\sqrt{N}}, \frac{1}{2\sqrt{N}}]^3} H(s,u)du\right)\left(\eta_{Ns-}(y)-\eta_{Ns-}(x)\right)d\hat{N}_{Ns}^{x,y}
\end{align*}
for any $H\in C_c\left([0, T]\times \mathbb{R}^3\right)$. We denote by $\mathcal{Y}$ the Gaussian time-space white noise on $[0, T]\times \mathbb{R}^3$ such that
$\mathcal{Y}(H)$ follows the normal distribution with mean zero and variance
\[
12p(1-p)\gamma_3\int_0^T\int_{\mathbb{R}^3}H^2(s,u)dsdu
\]
for any  $H\in C_c\left([0, T]\times \mathbb{R}^3\right)$.

For each $N\geq 1$ and $0\leq t\leq T$, we define
\[
b_t^N(s,u)=\sqrt{N}\sum_{x\in \mathbb{Z}^3}v\left(N(t-s), x\right)1_{\{u\in \frac{x}{\sqrt{N}}+(-\frac{1}{2\sqrt{N}}, \frac{1}{2\sqrt{N}}]^3, s\leq t\}}
\]
for any $(s, u)\in [0, T]\times \mathbb{R}^3$. By \eqref{equ 4.2}, we have
\begin{equation}\label{equ 4.4}
\frac{1}{N^{\frac{3}{4}}}\mathcal{M}_{tN}^{tN}=\mathcal{Y}^N(b_t^N)
\end{equation}
for any $0\leq t\leq T$. It is shown in \cite{Birkner2007}, via the local central limit theorem of the simple random walk, that
$b_t^N$ converges, as $N\rightarrow+\infty$, to $b_t$ uniformly and in $L^2$ on $[0, T]\times \mathbb{R}^3$, where
\[
b_t(s,u)=\int_0^{t-s}\frac{1}{(4\pi r)^{\frac{3}{2}}}e^{-\frac{\|u\|_2^2}{4r}}dr1_{\{s\leq t\}},
\]
where $\|u\|_2$ is the $l_2$-norm of $u$. Hence, roughly speaking, Lemma \ref{lemma 4.3} holds with $\Xi_j=\mathcal{Y}(b_{t_j})$ when we can show that $\mathcal{Y}^N$ converges to $\mathcal{Y}$. So we require the following lemma.

\begin{lemma}\label{lemma 4.4}
Let $\eta_0$ be distributed with $\nu_p$, then for any $H\in C_c\left([0, T]\times \mathbb{R}^3\right)$, $\mathcal{Y}^N(H)$ converges weakly to $\mathcal{Y}(H)$ as $N\rightarrow+\infty$ and furthermore,
\[
\lim_{N\rightarrow+\infty}\mathbb{E}_{\nu_p}\left(\left(\mathcal{Y}^N(H)\right)^2\right)=\mathbb{E}\left(\left(\mathcal{Y}(H)\right)^2\right).
\]
\end{lemma}

\proof[Proofs of Lemmas \ref{lemma 4.4} and \ref{lemma 4.3}]

We first prove Lemma \ref{lemma 4.4}. For any given integer $k\geq 1$ and $A_1, A_2,\ldots, A_k\subseteq \mathbb{R}^3$ which are disjoint bounded measurable domains, we let
\[
\varpi_{t,j}^N=\mathcal{Y}^N([0, t]\times A_j)
\]
for $1\leq j\leq k, 0\leq t\leq T$. We further define $\{\varpi_{t,j}:~0\leq t\leq T\}_{1\leq j\leq k}$ as independent Brownian motions such that
\[
{\rm Var}(\varpi_{t,j})=12p(1-p)\gamma_3\mathfrak{m}(A_j)t
\]
for all $1\leq j\leq k$, where $\mathfrak{m}(A_j)$ is the Lebesgue measure of $A_j$. According to the definition of $\mathcal{Y}$, to complete the proof of Lemma \ref{lemma 4.4}, we only need to show that, for any $1\leq j,l\leq k$, the cross variation process $\langle \varpi_j^N, \varpi_l^N \rangle_t$ of $\varpi_{t,j}^N$ and $\varpi_{t,l}^N$ converges in $L^2$ to
\[
12p(1-p)\gamma_3\mathfrak{m}(A_j)t1_{\{l=j\}}
\]
as $N\rightarrow+\infty$. By \eqref{equ 4.cross variation}, $\langle \varpi_j^N, \varpi_l^N \rangle_t=0$ when $j\neq l$ since $A_j$ and $A_l$ are disjoint. Hence we only deal with the $l=j$ case. According to the definition of $\mathcal{Y}^N$, we have
\[
\langle \varpi_j^N \rangle_t=N^{\frac{1}{2}}\sum_{x\in \mathbb{Z}^3}\sum_{y:y\sim x}\int_0^t\left(\int_{\frac{x}{\sqrt{N}}+(-\frac{1}{2\sqrt{N}}, \frac{1}{2\sqrt{N}}]^3}1_{\{u\in A_j\}}du\right)^2\left(\eta_{Ns}(y)-\eta_{Ns}(x)\right)^2 Nds.
\]
According to \eqref{equ duality minus square}, for $x\sim y$,
\[
\lim_{N\rightarrow+\infty}\mathbb{E}_{\nu_p}\left(\left(\eta_{Ns}(y)-\eta_{Ns}(x)\right)^2\right)
=2p(1-p)\hat{\mathbb{P}}(\tau_{xy}=+\infty)=2p(1-p)\gamma_3
\]
and hence
\begin{align}\label{equ 4.5 two}
\lim_{N\rightarrow+\infty}\mathbb{E}_{\nu_p}\langle \varpi_j^N \rangle_t&=2p(1-p)\gamma_3\times 6t \lim_{N\rightarrow+\infty}N^{\frac{3}{2}}
\sum_{x\in \mathbb{Z}^3}\left(\int_{\frac{x}{\sqrt{N}}+(-\frac{1}{2\sqrt{N}}, \frac{1}{2\sqrt{N}}]^3}1_{\{u\in A_j\}}du\right)^2 \notag\\
&=12p(1-p)\gamma_3\mathfrak{m}(A_j)t.
\end{align}
Therefore, to finish the proof of Lemma \ref{lemma 4.4}, we only need to show that
\begin{equation}\label{equ 4.6}
\lim_{N\rightarrow+\infty}{\rm Var}_{\nu_p}\left(\langle \varpi_j^N \rangle_t\right)=0.
\end{equation}
According to the expression of $\langle \varpi_j^N \rangle_t$ and the bilinear property of the covariance operator, it is not difficult to check that \eqref{equ 4.6} holds if
\begin{equation}\label{equ 4.7}
\lim_{r\rightarrow+\infty}\sup_{x\sim y, z\sim w, s\geq 0}{\rm Cov}\left(\left(\eta_{s}(y)-\eta_{s}(x)\right)^2, \left(\eta_{s+r}(z)-\eta_{s+r}(w)\right)^2\right)=0.
\end{equation}
Now we check \eqref{equ 4.7}. Let $(t_1, t_2, t_3, t_4)=(s,s,s+r, s+r)$, $(x_1, x_2, x_3, x_4)=(z,w,x,y)$ and
\[
\left(\mathcal{X}^{x_1, t_4-t_4}_{t,1}, \mathcal{X}^{x_2, t_4-t_3}_{t,2}, \mathcal{X}^{x_3, t_4-t_2}_{t,3}, \mathcal{X}^{x_4, t_4-t_1}_{t,4}\right)
=\left(\mathcal{X}^{z, 0}_{t,1}, \mathcal{X}^{w, 0}_{t,2}, \mathcal{X}^{x, r}_{t,3}, \mathcal{X}^{y, r}_{t,4}\right)
\]
defined as in Section \ref{section three}. We define $\beta$ as
\[
\beta=\inf\left\{t:~t>r, \left\{\mathcal{X}^{z, 0}_{t,1}, \mathcal{X}^{w, 0}_{t,2}\right\}\bigcap\left\{\mathcal{X}^{x, r}_{t,3}, \mathcal{X}^{y, r}_{t,4}\right\}\neq \emptyset\right\},
\]
then we denote by $\left\{\left(\hat{\mathcal{X}}^{x, r}_{t,3}, \hat{\mathcal{X}}^{y, r}_{t,4}\right)\right\}_{t\geq 0}$ the coalescing random walk starting from $(x, y)$ such that $\left\{\left(\hat{\mathcal{X}}^{x, r}_{t,3}, \hat{\mathcal{X}}^{y, r}_{t,4}\right)\right\}_{t\geq 0}$ is independent of $\left\{\left(\mathcal{X}^{z, 0}_{t,1}, \mathcal{X}^{w, 0}_{t,2}\right)\right\}_{t\geq 0}$ and
\[
\left(\hat{\mathcal{X}}^{x, r}_{t,3}, \hat{\mathcal{X}}^{y, r}_{t,4}\right)=\left(\mathcal{X}^{x, r}_{t,3}, \mathcal{X}^{y, r}_{t,4}\right)
\]
for $0\leq t\leq \beta$. Hence, $\left\{\left(\hat{\mathcal{X}}^{x, r}_{t,3}, \hat{\mathcal{X}}^{y, r}_{t,4}\right)\right\}_{t\geq 0}$ and $\left\{\left(\mathcal{X}^{x, r}_{t,3}, \mathcal{X}^{y, r}_{t,4}\right)\right\}_{t\geq 0}$ have the same distribution. Then, by Proposition \ref{proposition duality},
\begin{align}\label{equ 4.8}
&\mathbb{E}_{\nu_p}\left(\left(\eta_{s+r}(z)-\eta_{s+r}(w)\right)^2\right)\mathbb{E}_{\nu_p}\left(\left(\eta_{s}(y)-\eta_{s}(x)\right)^2\right)\notag\\
&=\hat{\mathbb{E}}\left(\int_{\{0, 1\}^{\mathbb{Z}^3}}\left(\eta\left(\mathcal{X}^{x, r}_{s+r,3}\right)-\eta\left(\mathcal{X}^{y, r}_{s+r,4}\right)\right)^2\nu_p(d\eta)\right)\notag\\
&\text{\quad\quad}\times \hat{\mathbb{E}}\left(\int_{\{0, 1\}^{\mathbb{Z}^3}}\left(\eta\left(\mathcal{X}^{z, 0}_{s+r,1}\right)-\eta\left(\mathcal{X}^{w, r}_{s+r,2}\right)\right)^2\nu_p(d\eta)\right)\notag\\
&=\hat{\mathbb{E}}\left(\int_{\{0, 1\}^{\mathbb{Z}^3}}\left(\eta\left(\hat{\mathcal{X}}^{x, r}_{s+r,3}\right)-\eta\left(\hat{\mathcal{X}}^{y, r}_{s+r,4}\right)\right)^2\nu_p(d\eta)\right)\notag\\
&\text{\quad\quad}\times \hat{\mathbb{E}}\left(\int_{\{0, 1\}^{\mathbb{Z}^3}}\left(\eta\left(\mathcal{X}^{z, 0}_{s+r,1}\right)-\eta\left(\mathcal{X}^{w, r}_{s+r,2}\right)\right)^2\nu_p(d\eta)\right)\notag\\
&=\hat{\mathbb{E}}\Bigg(\int_{\{0, 1\}^{\mathbb{Z}^3}}\left(\eta\left(\hat{\mathcal{X}}^{x, r}_{s+r,3}\right)-\eta\left(\hat{\mathcal{X}}^{y, r}_{s+r,4}\right)\right)^2\nu_p(d\eta)\\
&\text{\quad\quad\quad}\times \int_{\{0, 1\}^{\mathbb{Z}^3}}\left(\eta\left(\mathcal{X}^{z, 0}_{s+r,1}\right)-\eta\left(\mathcal{X}^{w, r}_{s+r,2}\right)\right)^2\nu_p(d\eta)\Bigg)\notag
\end{align}
and
\begin{align*}
&\mathbb{E}_{\nu_p}\left(\left(\eta_{s+r}(z)-\eta_{s+r}(w)\right)^2\left(\eta_{s}(y)-\eta_{s}(x)\right)^2\right)\\
&=\hat{\mathbb{E}}\Bigg(\int_{\{0, 1\}^{\mathbb{Z}^3}}\left(\eta\left(\mathcal{X}^{x, r}_{s+r,3}\right)-\eta\left(\mathcal{X}^{y, r}_{s+r,4}\right)\right)^2\left(\eta\left(\mathcal{X}^{z, 0}_{s+r,1}\right)-\eta\left(\mathcal{X}^{w, r}_{s+r,2}\right)\right)^2\nu_p(d\eta)\Bigg).
\end{align*}
On the event $\{\beta>s+r\}$, we have
\begin{align*}
&\int_{\{0, 1\}^{\mathbb{Z}^3}}\left(\eta\left(\hat{\mathcal{X}}^{x, r}_{s+r,3}\right)-\eta\left(\hat{\mathcal{X}}^{y, r}_{s+r,4}\right)\right)^2\nu_p(d\eta)
\int_{\{0, 1\}^{\mathbb{Z}^3}}\left(\eta\left(\mathcal{X}^{z, 0}_{s+r,1}\right)-\eta\left(\mathcal{X}^{w, r}_{s+r,2}\right)\right)^2\nu_p(d\eta)\\
&=\int_{\{0, 1\}^{\mathbb{Z}^3}}\left(\eta\left(\hat{\mathcal{X}}^{x, r}_{s+r,3}\right)-\eta\left(\hat{\mathcal{X}}^{y, r}_{s+r,4}\right)\right)^2\left(\eta\left(\mathcal{X}^{z, 0}_{s+r,1}\right)-\eta\left(\mathcal{X}^{w, r}_{s+r,2}\right)\right)^2\nu_p(d\eta)\\
&=\int_{\{0, 1\}^{\mathbb{Z}^3}}\left(\eta\left(\mathcal{X}^{x, r}_{s+r,3}\right)-\eta\left(\mathcal{X}^{y, r}_{s+r,4}\right)\right)^2\left(\eta\left(\mathcal{X}^{z, 0}_{s+r,1}\right)-\eta\left(\mathcal{X}^{w, r}_{s+r,2}\right)\right)^2\nu_p(d\eta)
\end{align*}
since
\[
\left\{\mathcal{X}^{z, 0}_{t,1}, \mathcal{X}^{w, 0}_{t,2}\right\}\bigcap\left\{\hat{\mathcal{X}}^{x, r}_{t,3}, \hat{\mathcal{X}}^{y, r}_{t,4}\right\}=\emptyset.
\]
Therefore, by \eqref{equ 4.8},
\begin{equation}\label{equ 4.9}
\left|{\rm Cov}\left(\left(\eta_{s}(y)-\eta_{s}(x)\right)^2, \left(\eta_{s+r}(z)-\eta_{s+r}(w)\right)^2\right)\right|\leq 2\hat{\mathbb{P}}\left(\beta\leq s+r\right).
\end{equation}
According to the definition of $\beta$, we have $\beta=\inf\{\beta_{13}, \beta_{14}, \beta_{23}, \beta_{24}\}$, where
\[
\beta_{ij}=\inf\left\{t:~t>r, \mathcal{X}^{x_i, 0}_{t,i}=\mathcal{X}^{x_j, r}_{t,j}\right\}
\]
for $i\in \{1,2\}$ and $j\in \{3,4\}$. Note that $\mathcal{X}^{x, r}_{t,3}$ is frozen in $x$ for $0\leq t\leq r$, then $\mathcal{X}^{x, r}_{r+t,3}-\mathcal{X}^{z, 0}_{r+t,1}$ evolves as $X_{2t}$ until hits $O$. Therefore,
\[
\hat{\mathbb{P}}(\beta_{13}\leq s+r)\leq \sum_{v\in \mathbb{Z}^3}p_r(z,x+v)\hat{\mathbb{P}}\left(X_t=O\text{~for some~}t\geq 0\Big|X_0=v\right).
\]
According to the well-known result that $p_r(z,x+v)\leq p_r(O,O)$, we have
\[
\hat{\mathbb{P}}(\beta_{13}\leq s+r)\leq \sum_{v:\|v\|_1\leq M}p_r(O, O)+\sup_{v:\|v\|_1>M}\hat{\mathbb{P}}\left(X_t=O\text{~for some~}t\geq 0\Big|X_0=v\right)
\]
for any $M\geq 0$. According to similar analysis, the above inequality still holds when $\beta_{13}$ is replaced by $\beta_{23}, \beta_{14}, \beta_{24}$. Therefore, by \eqref{equ 4.9}
\begin{align}\label{equ 4.10}
&\limsup_{r\rightarrow+\infty}\sup_{x\sim y, z\sim w, s\geq 0}{\rm Cov}\left(\left(\eta_{s}(y)-\eta_{s}(x)\right)^2, \left(\eta_{s+r}(z)-\eta_{s+r}(w)\right)^2\right)\notag\\
&\leq 8\left(\sum_{v:\|v\|_1\leq M}\limsup_{r\rightarrow+\infty}p_r(O, O)+\sup_{v:\|v\|_1>M}\hat{\mathbb{P}}\left(X_t=O\text{~for some~}t\geq 0\Big|X_0=v\right)\right)\notag\\
&=8\sup_{v:\|v\|_1>M}\hat{\mathbb{P}}\left(X_t=O\text{~for some~}t\geq 0\Big|X_0=v\right).
\end{align}
Since the simple random walk on $\mathbb{Z}^3$ is transient, we have
\[
\lim_{M\rightarrow+\infty}\sup_{v:\|v\|_1>M}\hat{\mathbb{P}}\left(X_t=O\text{~for some~}t\geq 0\Big|X_0=v\right)=0
\]
and hence \eqref{equ 4.7} follows from \eqref{equ 4.10}. As we have explained, \eqref{equ 4.6} follows from \eqref{equ 4.7} and Lemma \ref{lemma 4.4} follows from \eqref{equ 4.6}. Hence, the proof of Lemma \ref{lemma 4.4} is complete.

\quad

As we have explained, $b_t^N$ converges to $b_t$ uniformly and in $L_2$ on $[0, T]\times \mathbb{R}^3$. Therefore, by Lemma \ref{lemma 4.4} and \eqref{equ 4.4}, Lemma \ref{lemma 4.3} holds with $\Xi_j=\mathcal{Y}(b_{t_j})$ for $1\leq j\leq m$. Hence, the proof of Lemma \ref{lemma 4.3} is complete.

\qed

The following lemma shows that any finite dimensional distribution of
\[
\left\{\frac{1}{N^{\frac{3}{4}}}\int^{tN}_0\left(\eta_s(O)-p\right)ds\right\}_{0\leq t\leq T}
\]
converges weakly to the corresponding finite dimensional distribution of $\{C_3\zeta_t\}_{0\leq t\leq T}$.

\begin{lemma}\label{lemma 4.5}
Let $\eta_0$ be distributed with $\nu_p$, then for any integer $m\geq 1$ and $0\leq t_1<t_2<\ldots<t_m\leq T$, $\frac{1}{N^{\frac{3}{4}}}\left(\int^{t_1N}_0\left(\eta_s(O)-p\right)ds, \ldots, \int^{t_mN}_0\left(\eta_s(O)-p\right)ds\right)$ converges weakly, as $N\rightarrow+\infty$, to $C_3\left(\zeta_{t_1},\cdots, \zeta_{t_m}\right)$.
\end{lemma}

\proof

By \eqref{equ 4.3} and Lemmas \ref{lemma small error}, \ref{lemma 4.3}, we only need to show that, for any $0\leq s<t\leq T$,
\begin{equation}\label{equ 4.11}
\lim_{N\rightarrow+\infty}{\rm Cov}\left(\frac{1}{N^{\frac{3}{4}}}\mathcal{M}_{sN}^{sN}, \frac{1}{N^{\frac{3}{4}}}\mathcal{M}_{tN}^{tN}\right)=C_3^2\left(s^{\frac{3}{2}}+t^{\frac{3}{2}}-\frac{1}{2}(t-s)^{\frac{3}{2}}-\frac{1}{2}(t+s)^{\frac{3}{2}}\right).
\end{equation}
Now we check \eqref{equ 4.11}. By \eqref{equ 4.2} and \eqref{equ duality minus square},
\begin{align*}
&{\rm Cov}\left(\frac{1}{N^{\frac{3}{4}}}\mathcal{M}_{sN}^{sN}, \frac{1}{N^{\frac{3}{4}}}\mathcal{M}_{tN}^{tN}\right)\\
&=\frac{1}{N^{\frac{3}{2}}}\sum_{x\in \mathbb{Z}^3}\sum_{y:y\sim x}\int_0^{sN}v(tN-r, x)v(sN-r, x)\mathbb{E}_{\nu_p}\left((\eta_r(y)-\eta_r(x))^2\right)dr \\
&=\frac{2p(1-p)}{N^{\frac{1}{2}}}\sum_{x\in \mathbb{Z}^3}\sum_{y:y\sim x}\int_0^sv((t-\theta)N, x)v((s-\theta)N, x)\hat{\mathbb{P}}\left(\tau_{xy}>\theta N\right)d\theta.
\end{align*}
Since $\lim_{N\rightarrow+\infty}\hat{\mathbb{P}}\left(\tau_{xy}>\theta N\right)=\gamma_3$ for $x\sim y$, to check \eqref{equ 4.11}, we only need to show that
\begin{align}\label{equ 4.12}
&\lim_{N\rightarrow+\infty}\frac{1}{N^{\frac{1}{2}}}\sum_{x\in \mathbb{Z}^3}\int_0^sv((t-\theta)N, x)v((s-\theta)N, x)d\theta\\
&=\frac{1}{3\pi^{\frac{3}{2}}}\left(s^{\frac{3}{2}}+t^{\frac{3}{2}}-\frac{1}{2}(t-s)^{\frac{3}{2}}-\frac{1}{2}(t+s)^{\frac{3}{2}}\right).\notag
\end{align}
According to the definition of $v(t,x)$, we have
\[
\sum_{x\in \mathbb{Z}^3}\int_0^sv((t-\theta)N, x)v((s-\theta)N, x)d\theta
=\int_0^s\left(\int_0^{(t-\theta)N}\int_0^{(s-\theta)N}p_{r_1+r_2}(O, O)dr_1dr_2\right)d\theta.
\]
According to the local central limit theorem of the simple random walk on $\mathbb{Z}^d$, we have
\[
\lim_{t\rightarrow+\infty}\frac{p_t(O, O)}{(4\pi t)^{-\frac{3}{2}}}=1.
\]
Hence, for any $\epsilon>0$, there exists $M>0$ such that
\[
(4\pi t)^{-\frac{3}{2}}(1-\epsilon)\leq p_t(O, O)\leq (4\pi t)^{-\frac{3}{2}}(1+\epsilon)
\]
when $t\geq M$. Therefore, for sufficiently large $N$,
\begin{equation}\label{equ 4.13}
K_1^N\leq \int_0^s\left(\int_0^{(t-\theta)N}\int_0^{(s-\theta)N}p_{r_1+r_2}(O, O)dr_1dr_2\right)d\theta\leq K_2^N,
\end{equation}
where
\[
K_1^N=(1-\epsilon)\int_0^s\left(\int_M^{(t-\theta)N}\int_0^{(s-\theta)N}\frac{1}{(r_1+r_2)^{\frac{3}{2}}8\pi^{\frac{3}{2}}}dr_1dr_2\right)d\theta
\]
and
\[
K_2^N=sM\int_0^{+\infty}p_\theta(O,O)d\theta+(1+\epsilon)\int_0^s\left(\int_M^{(t-\theta)N}\int_0^{(s-\theta)N}\frac{1}{(r_1+r_2)^{\frac{3}{2}}8\pi^{\frac{3}{2}}}
dr_1dr_2\right)d\theta.
\]
Note that $\int_0^{+\infty}p_\theta(O,O)d\theta<+\infty$ since the simple random walk on $\mathbb{Z}^3$ is transient. By direct calculation, it is not difficult to show that
\[
\lim_{N\rightarrow+\infty}\frac{K_1^N}{\sqrt{N}}=\frac{1}{3\pi^{\frac{2}{3}}}\left(s^{\frac{3}{2}}+t^{\frac{3}{2}}-\frac{1}{2}(t-s)^{\frac{3}{2}}-\frac{1}{2}(t+s)^{\frac{3}{2}}\right)(1-\epsilon)
\]
and
\[
\lim_{N\rightarrow+\infty}\frac{K_2^N}{\sqrt{N}}=\frac{1}{3\pi^{\frac{2}{3}}}\left(s^{\frac{3}{2}}+t^{\frac{3}{2}}-\frac{1}{2}(t-s)^{\frac{3}{2}}-\frac{1}{2}(t+s)^{\frac{3}{2}}\right)(1+\epsilon).
\]
Since $\epsilon$ is arbitrary, \eqref{equ 4.12} follows from \eqref{equ 4.13} and the proof is complete.

\qed

At last, we prove Theorem \ref{theorem 2.1 main sample path CLT} in case $d=3$.

\proof[Proof of Theorem \ref{theorem 2.1 main sample path CLT} in case $d=3$]

By Lemma \ref{lemma 4.5}, we only need to show that
\[
\left\{\frac{1}{N^{\frac{3}{4}}}\int_0^{tN}(\eta_s(O)-p)ds:~0\leq t\leq T\right\}_{N\geq 1}
\]
are tight under the Skorohod topology. A well known criterion of this tightness is to show that there exist $a>0, b>0$ and $K>0$ independent of $s, t>0$ such that
\begin{equation}\label{equ criterion of tightness}
\mathbb{E}_{\nu_p}\left(\left|\frac{1}{N^{\frac{3}{4}}}\int_{sN}^{tN}(\eta_r(O)-p)dr\right|^a\right)\leq K|t-s|^{1+b}
\end{equation}
for any $t,s\geq 0$. Here we show that \eqref{equ criterion of tightness} holds with $a=2$ and $b=\frac{1}{2}$. According to \eqref{equ conservative mean}, for $t>s$, we have
\[
\mathbb{E}_{\nu_p}\left(\left|\frac{1}{N^{\frac{3}{4}}}\int_{sN}^{tN}(\eta_r(O)-p)dr\right|^2\right)
=\frac{2}{N^{\frac{3}{2}}}\int_{sN}^{tN}\left(\int_{sN}^{\theta}{\rm Cov}_{\nu_p}\left(\eta_r(O), \eta_\theta(O)\right)dr\right)d\theta.
\]
According to \eqref{equ dual eta(x)}, \eqref{equ conservative mean} and the strong Markov property, we have
\begin{align*}
{\rm Cov}_{\nu_p}\left(\eta_r(O), \eta_\theta(O)\right)&=\mathbb{E}_{\nu_p}\left(\eta_r(O)\eta_\theta(O)\right)-p^2\\
&=\mathbb{E}_{\nu_p}\left(\eta_r(O)\mathbb{E}_{\eta_r}\eta_{\theta-r}(O)\right)-p^2\\
&=\sum_{x\in \mathbb{Z}^3}p_{\theta-r}(O, x)\mathbb{E}_{\nu_p}\left(\eta_r(O)\eta_r(x)\right)-p^2\\
&=\sum_{x\in \mathbb{Z}^3}p_{\theta-r}(O, x)\left(\mathbb{E}_{\nu_p}\left(\eta_r(O)\eta_r(x)\right)-p^2\right).
\end{align*}
By \eqref{equ duality relationship},
\[
\mathbb{E}_{\nu_p}\left(\eta_r(O)\eta_r(x)\right)=ph(r,x)+p^2(1-h(r,x)),
\]
where $h(r,x)=\hat{\mathbb{P}}\left(X_t=O\text{~for some~}t<r\Big|X_0=x\right)$. Hence,
\begin{align*}
\mathbb{E}_{\nu_p}\left(\eta_r(O)\eta_r(x)\right)-p^2&=p(1-p)h(r,x)\leq p(1-p)h(+\infty, x)\\
&=p(1-p)\frac{\int_0^{+\infty}p_v(x, O)dv}{\int_0^{+\infty}p_v(O,O)dv}=6p(1-p)\gamma_3\int_0^{+\infty}p_v(x, O)dv.
\end{align*}
In conclusion, by local central limit theorems of simple random walks on $\mathbb{Z}^3$,
\begin{align*}
&\mathbb{E}_{\nu_p}\left(\left|\frac{1}{N^{\frac{3}{4}}}\int_{sN}^{tN}(\eta_r(O)-p)dr\right|^2\right)\\
&\leq \frac{12p(1-p)\gamma_3}{N^{\frac{3}{2}}}\int_{sN}^{tN}\left(\int_{sN}^\theta\sum_{x\in \mathbb{Z}^3}p_{\theta-r}(O,x)\left(\int_0^{+\infty}p_v(x, O)dv\right)dr\right)d\theta\\
&=\frac{12p(1-p)\gamma_3}{N^{\frac{3}{2}}}\int_{sN}^{tN}\left(\int_{sN}^\theta\left(\int_0^{+\infty}p_{\theta-r+v}(O, O)dv\right)dr\right)d\theta\\
&\leq \frac{12p(1-p)\gamma_3K_2}{N^{\frac{3}{2}}}\int_{sN}^{tN}\left(\int_{sN}^\theta\left(\int_0^{+\infty}\left(\theta-r+v\right)^{-\frac{3}{2}}dv\right)dr\right)
d\theta\\
&=\frac{12p(1-p)\gamma_3K_2}{N^{\frac{3}{2}}}\frac{8}{3}(t-s)^{\frac{3}{2}}N^{\frac{3}{2}}=32p(1-p)\gamma_3K_2(t-s)^{\frac{3}{2}}
\end{align*}
for some $K_2<+\infty$ independent of $s, t$. Hence, Equation \eqref{equ criterion of tightness} holds for $a=2, b=\frac{1}{2}$ and the proof is complete.

\qed

\section{The proof of Theorem \ref{theorem 2.1 main sample path CLT}: $d=4$ case}\label{section six d=4}
In this section, we prove Theorem \ref{theorem 2.1 main sample path CLT} in the $d=4$ case. We first introduce some notations and definitions for later use. For each $N\geq 1$ and $x\in \mathbb{Z}^4$, we define
\[
\phi_N(x)=\int_0^{+\infty}e^{-\frac{s}{N}}p_s(O, x)ds.
\]
Then for all $N\geq 1$ and $\eta\in \{0, 1\}^{\mathbb{Z}^4}$, we define
\[
\mathcal{H}_N(\eta)=\sum_{x\in \mathbb{Z}^4}\left(\eta(x)-p\right)\phi_N(x).
\]
According to an analysis similar with that leading to \eqref{equ generator and occupation time}, we have
\begin{equation}\label{equ d=4 1}
\mathcal{L}\mathcal{H}_N(\eta)=\frac{1}{N}\mathcal{H}_N(\eta)-(\eta(O)-p).
\end{equation}
For any $t\geq 0$, we define
\[
\widetilde{\mathcal{M}}_t^N=\mathcal{H}_N(\eta_t)-\mathcal{H}_N(\eta_0)-\int_0^t\mathcal{L}\mathcal{H}_N(\eta_s)ds,
\]
then, according to the Dynkin's martingale formula, $\{\widetilde{\mathcal{M}}_t^N\}_{t\geq 0}$ is a martingale. By \eqref{equ d=4 1}, we have

\begin{equation}\label{equ d=4 2}
\frac{\int_0^{tN}\left(\eta_s(O)-p\right)ds}{\sqrt{N\log N}}
=\frac{\widetilde{\mathcal{M}}_{tN}^N}{\sqrt{N\log N}}-\frac{\mathcal{H}_N(\eta_{tN})}{\sqrt{N\log N}}+\frac{\mathcal{H}_N(\eta_{0})}{\sqrt{N\log N}}+\frac{\int_0^{tN}\mathcal{H}_N(\eta_s)ds}{\sqrt{N\log N}N}.
\end{equation}
The following lemma shows that the term $-\frac{\mathcal{H}_N(\eta_{tN})}{\sqrt{N\log N}}+\frac{\mathcal{H}_N(\eta_{0})}{\sqrt{N\log N}}+\frac{\int_0^{tN}\mathcal{H}_N(\eta_s)ds}{\sqrt{N\log N}N}$ is a small error.

\begin{lemma}\label{lemma d=4 1 small error}
Let $\eta_0$ be distributed with $\nu_p$. For any given $t\geq 0$,
\[
\lim_{N\rightarrow+\infty}\left(-\frac{\mathcal{H}_N(\eta_{tN})}{\sqrt{N\log N}}+\frac{\mathcal{H}_N(\eta_{0})}{\sqrt{N\log N}}+\frac{\int_0^{tN}\mathcal{H}_N(\eta_s)ds}{\sqrt{N\log N}N}\right)=0
\]
in $L^2$.
\end{lemma}

\proof

According to the Cauchy-Schwartz inequality, we only need to show that
\begin{equation}\label{equ d=4 3}
\lim_{N\rightarrow+\infty}\frac{1}{N\log N}\left(\sup_{s\geq 0}\mathbb{E}_{\nu_p}\left(\left(\mathcal{H}_N(\eta_s)\right)^2\right)\right)=0.
\end{equation}

By \eqref{equ conservative mean}, \eqref{equ voter covariance} and the strong Markov property of the simple random walk, we have
\begin{align*}
\mathbb{E}_{\nu_p}\left(\left(\mathcal{H}_N(\eta_s)\right)^2\right)&={\rm Cov}_{\nu_p}\left(\mathcal{H}_N(\eta_s), \mathcal{H}_N(\eta_s)\right)\\
&=\sum_{x\in \mathbb{Z}^4}\sum_{y\in \mathbb{Z}^4}\phi_N(x)\phi_N(y){\rm Cov}_{\nu_p}\left(\eta_s(x), \eta_s(y)\right)\\
&=p(1-p)\sum_{x\in \mathbb{Z}^4}\sum_{y\in \mathbb{Z}^4}\phi_N(x)\phi_N(y)\hat{\mathbb{P}}(\tau_{xy}\leq s)\\
&\leq p(1-p)\sum_{x\in \mathbb{Z}^4}\sum_{y\in \mathbb{Z}^4}\phi_N(x)\phi_N(y)\hat{\mathbb{P}}(\tau_{xy}<+\infty)\\
&=p(1-p)\sum_{x\in \mathbb{Z}^4}\sum_{y\in \mathbb{Z}^4}\phi_N(x)\phi_N(y)\frac{\int_0^{+\infty} p_r(x,y)dr}{\int_0^{+\infty}p_r(O, O)dr}\\
&=8p(1-p)\gamma_4\sum_{x\in \mathbb{Z}^4}\sum_{y\in \mathbb{Z}^4}\phi_N(x)\phi_N(y)\int_0^{+\infty} p_r(x,y)dr\\
&=8p(1-p)\gamma_4\int_0^{+\infty}\int_0^{+\infty}e^{-\frac{s_1+s_2}{N}}\left(\int_0^{+\infty}p_{r+s_1+s_2}(O, O)dr\right)ds_1ds_2\\
&=8p(1-p)\gamma_4\int_0^{+\infty}N^2\theta e^{-\theta}\left(\int_0^{+\infty}p_{r+N\theta}(O, O)dr\right)d\theta.
\end{align*}
According to the local central limit theorem of the simple random walk, there exists $K_4<+\infty$ independent of $\theta, r, N$ such that
\[
p_{r+N\theta}(O, O)\leq \frac{K_4}{(r+N\theta)^2}.
\]
Therefore,
\begin{align*}
\mathbb{E}_{\nu_p}\left(\left(\mathcal{H}_N(\eta_s)\right)^2\right)&\leq 8p(1-p)\gamma_4 K_4 \int_0^{+\infty}N^2\theta e^{-\theta}\frac{1}{N\theta}d\theta\\
&=8Np(1-p)\gamma_4K_4\int_0^{+\infty}e^{-\theta} d\theta=8Np(1-p)\gamma_4K_4
\end{align*}
and hence
\[
\frac{1}{N\log N}\left(\sup_{s\geq 0}\mathbb{E}_{\nu_p}\left(\left(\mathcal{H}_N(\eta_s)\right)^2\right)\right)\leq \frac{8p(1-p)\gamma_4K_4}{\log N}.
\]
As a result, \eqref{equ d=4 3} holds and the proof is complete.

\qed

The following lemma shows that any finite dimensional distribution of $\left\{\frac{\int_0^{tN}\left(\eta_s(O)-p\right)ds}{h_4(N)}\right\}_{0\leq t\leq T}$ converges weakly to the corresponding finite dimensional distribution of $\left\{C_4W_t\right\}_{0\leq t\leq T}$.

\begin{lemma}\label{lemma d=4 finite dimension}
Let $\eta_0$ be distributed with $\nu_p$, then for each integer $m\geq 1$ and $0\leq t_1<t_2<\ldots<t_m\leq T$,
\[
\frac{1}{\sqrt{N\log N}}\left(\int_0^{t_1N}\left(\eta_s(O)-p\right)ds, \int_0^{t_2N}\left(\eta_s(O)-p\right)ds,\ldots, \int_0^{t_mN}\left(\eta_s(O)-p\right)ds\right)
\]
converges weakly, as $N\rightarrow+\infty$, to $C_4\left(W_{t_1}, W_{t_2}, \ldots, W_{t_m}\right)$.
\end{lemma}

As a preliminary of the proof of Lemma \ref{lemma d=4 finite dimension}, we require the following lemma.

\begin{lemma}\label{lemma d=4 convergence of mathcalH} We have
\[
\lim_{N\rightarrow+\infty}\frac{1}{\log N}\sum_{x\in \mathbb{Z}^4}\phi_N^2(x)=\frac{1}{16\pi^2}.
\]
\end{lemma}

\proof

According to the definition of $\phi_N$, we have
\begin{align*}
\sum_{x\in \mathbb{Z}^4}\phi_N^2(x)&=\int_0^{+\infty}\int_0^{+\infty}e^{-\frac{s_1+s_2}{N}}p_{s_1+s_2}(O, O)ds_1ds_2\\
&=\int_0^{+\infty}se^{-\frac{s}{N}}p_s(O, O)ds=N^2\int_0^{+\infty}se^{-s}p_{sN}(O, O)ds.
\end{align*}

According to the local central limit theorem of the simple random walk, for any $\epsilon>0$, there exists $M<+\infty$ depending on $\epsilon$ such that
\[
\frac{1}{16\pi^2}(1-\epsilon)t^{-2}\leq p_t(O,O)\leq \frac{1}{16\pi^2}(1+\epsilon)t^{-2}
\]
for all $t\geq M$. Hence, for sufficiently large $N$,
\begin{align*}
\sum_{x\in \mathbb{Z}^4}\phi_N^2(x)&\geq N^2\int_{\frac{M}{N}}^{+\infty}se^{-s}p_{sN}(O, O)ds\\
&\geq N^2\frac{1}{16\pi^2}(1-\epsilon)\int_{\frac{M}{N}}^{+\infty}se^{-s}\frac{1}{s^2N^2}ds\\
&=\frac{1}{16\pi^2}(1-\epsilon)\int_{\frac{M}{N}}^{+\infty}s^{-1}e^{-s}ds\\
&\geq \frac{1}{16\pi^2}(1-\epsilon)\int_{\frac{M}{N}}^\epsilon s^{-1}e^{-s}ds\geq \frac{1}{16\pi^2}(1-\epsilon)e^{-\epsilon}\int_{\frac{M}{N}}^\epsilon s^{-1}ds\\
&=\frac{1}{16\pi^2}(1-\epsilon)e^{-\epsilon}(\log \epsilon-\log M+\log N).
\end{align*}
Therefore, $\liminf_{N\rightarrow+\infty}\frac{1}{\log N}\sum_{x\in \mathbb{Z}^4}\phi_N^2(x)\geq \frac{1}{16\pi^2}(1-\epsilon)e^{-\epsilon}$. Since $\epsilon$ is arbitrary, let $\epsilon\rightarrow 0$ and then
\begin{equation}\label{equ lemma 6.3 liminf}
\liminf_{N\rightarrow+\infty}\frac{1}{\log N}\sum_{x\in \mathbb{Z}^4}\phi_N^2(x)\geq \frac{1}{16\pi^2}.
\end{equation}
Similarly, for sufficiently large $N$, we have
\begin{align*}
\sum_{x\in \mathbb{Z}^4}\phi_N^2(x)&\leq N^2\int_{\frac{M}{N}}^{+\infty}se^{-s}p_{sN}(O, O)ds+N^2\int_0^{\frac{M}{N}}sds\\
&\leq N^2\frac{1}{16\pi^2}(1+\epsilon)\int_{\frac{M}{N}}^{+\infty}se^{-s}\frac{1}{s^2N^2}ds+\frac{M^2}{2}\\
&=\frac{1}{16\pi^2}(1+\epsilon)\int_{\frac{M}{N}}^{+\infty}s^{-1}e^{-s}ds+\frac{M^2}{2}\\
&\leq \frac{1+\epsilon}{16\pi^2}\int_{\frac{M}{N}}^\epsilon s^{-1}ds+\frac{1+\epsilon}{16\pi^2\epsilon}\int_\epsilon^{+\infty}e^{-s}ds+\frac{M^2}{2}\\
&\leq \frac{1+\epsilon}{16\pi^2}\left(\log \epsilon-\log M+\log N\right)+\frac{1+\epsilon}{16\pi^2\epsilon}+\frac{M^2}{2}.
\end{align*}
Hence, $\limsup_{N\rightarrow+\infty}\frac{1}{\log N}\sum_{x\in \mathbb{Z}^4}\phi_N^2(x)\leq \frac{1+\epsilon}{16\pi^2}$. Since $\epsilon$ is arbitrary, let $\epsilon\rightarrow 0$ and then
\begin{equation}\label{equ lemma 6.3 limsup}
\limsup_{N\rightarrow+\infty}\frac{1}{\log N}\sum_{x\in \mathbb{Z}^4}\phi_N^2(x)\leq \frac{1}{16\pi^2}.
\end{equation}
Lemma \ref{lemma d=4 convergence of mathcalH} follows from \eqref{equ lemma 6.3 liminf} and \eqref{equ lemma 6.3 limsup}.

\qed

Now we give the proof of Lemma \ref{lemma d=4 finite dimension}.

\proof[The proof of Lemma \ref{lemma d=4 finite dimension}]

We denote by $\{\langle \widetilde{\mathcal{M}}^N\rangle_t\}_{t\geq 0}$ the quadratic variation process of the martingale $\{\widetilde{\mathcal{M}}^N_t\}_{t\geq 0}$. By \eqref{equ d=4 2} and Lemma \ref{lemma d=4 1 small error}, to prove Lemma \ref{lemma d=4 finite dimension}, we only need to show that
\begin{equation}\label{equ d=4 convergence of quadratic}
\lim_{N\rightarrow+\infty}\frac{1}{N\log N}\langle \widetilde{\mathcal{M}}^N\rangle_{tN}=C_4^2 t
\end{equation}
in $L^2$ for any $t\geq 0$. By Dynkin's martingale formula, we have
\begin{align*}
\frac{1}{N\log N}\langle \widetilde{\mathcal{M}}^N\rangle_{tN}&=\frac{1}{N\log N}\int_0^{tN}\left(\mathcal{L}\mathcal{H}_N^2(\eta_s)-2\mathcal{H}_N(\eta_s)\mathcal{L}\mathcal{H}_N(\eta_s)\right)ds\\
&=\frac{1}{N\log N}\sum_{x\in \mathbb{Z}^4}\sum_{y\sim x}\int_0^{tN}\phi_N^2(x)\left(\eta_{s}(y)-\eta_{s}(x)\right)^2ds\\
&=\frac{1}{\log N}\sum_{x\in \mathbb{Z}^4}\sum_{y\sim x}\int_0^{t}\phi_N^2(x)\left(\eta_{sN}(y)-\eta_{sN}(x)\right)^2ds.
\end{align*}
Then, by \eqref{equ duality minus square} and Lemma \ref{lemma d=4 convergence of mathcalH},
\begin{align*}
\lim_{N\rightarrow+\infty}\mathbb{E}_{\nu_p}\left(\frac{1}{N\log N}\langle \widetilde{\mathcal{M}}^N\rangle_{tN}\right)
&=8\times 2p(1-p)\hat{\mathbb{P}}\left(X_t\neq O\text{~for all~}t\Big|\|X_0\|=1\right)\frac{t}{16\pi^2}\\
&=\frac{16p(1-p)\gamma_4t}{16\pi^2}=C_4^2t.
\end{align*}
Hence, to prove \eqref{equ d=4 convergence of quadratic}, we only need to show that
\begin{equation}\label{equ 6.7}
\lim_{N\rightarrow+\infty}\frac{1}{N^2(\log N)^2}{\rm Var}_{\nu_p}\left(\langle \widetilde{\mathcal{M}}^N\rangle_{tN}\right)=0.
\end{equation}

The proof of \eqref{equ 6.7} is similar with that of \eqref{equ 4.6} and hence here we only give an outline. By utilizing the bilinear property of the covariance operator, it is easy to check that \eqref{equ 6.7} follows from a $d=4$ version of \eqref{equ 4.7}. The proof of \eqref{equ 4.7} given in Section \ref{section four d=3} relies on the transience of the simple random walk on the lattice, which can be extended to all $d\geq 4$ cases. Hence, \eqref{equ 6.7} holds and the proof of Lemma \ref{lemma d=4 finite dimension} is complete.

\qed

At last, we prove Theorem \ref{theorem 2.1 main sample path CLT} in the $d=4$ case.

\proof[Proof of Theorem \ref{theorem 2.1 main sample path CLT} in case $d=4$]

By Lemma \ref{lemma d=4 finite dimension}, we only need to show that
\[
\left\{\frac{1}{\sqrt{N\log N}}\int_0^{tN}\left(\eta_s(O)-p\right)ds:~0\leq t\leq T\right\}_{N\geq 1}
\]
are tight under the Skorohod topology. As in the proof of the $d=3$ case, we will show that there exists $K_5<+\infty$ independent of $N, s, t$ such that
\begin{equation}\label{equ 6.8}
\mathbb{E}_{\nu_p}\left(\left(\frac{1}{\sqrt{N\log N}}\int_{sN}^{tN}\left(\eta_r(O)-p\right)dr\right)^4\right)\leq K_5(t-s)^2
\end{equation}
for all $N\geq \max\{4, T\}$ and $0\leq s<t\leq T$ to complete the check of the aforesaid tightness. By direct calculation,
\begin{align*}
&\mathbb{E}_{\nu_p}\left(\left(\frac{1}{\sqrt{N\log N}}\int_{sN}^{tN}\left(\eta_r(O)-p\right)dr\right)^4\right)\\
&=\frac{1}{N^2(\log N)^2}\int_{sN}^{tN}\int_{sN}^{tN}\int_{sN}^{tN}\int_{sN}^{tN}\mathbb{E}_{\nu_p}\left(\prod_{j=1}^4\left(\eta_{t_j}(O)-p\right)\right)
dt_1dt_2dt_3dt_4.
\end{align*}
Then, by \eqref{equ 3.6} and \eqref{equ 3.7}, to prove \eqref{equ 6.8} we only need to show that there exists $K_6<+\infty$ independent of $N, s, t$ such that
\begin{equation}\label{equ 6.9}
\frac{1}{N\log N}\int_{sN}^{tN}\left(\int_{sN}^r\hat{\mathbb{P}}\left(\hat{\tau}_{34}(r, \theta)\leq r\right)d\theta\right)dr\leq K_6(t-s)
\end{equation}
for all $N\geq \max\{4, T\}$ and $0\leq s<t\leq T$, where $\hat{\tau}_{34}(r, \theta)$ is the $\hat{\tau}_{34}$ defined in Section \ref{section three} with $t_3=\theta, t_4=r$ and $(x_3, x_4)=(O, O)$.

According to the strong Markov property of the coalescing random walk,
\begin{align*}
\hat{\mathbb{P}}\left(\hat{\tau}_{34}(r, \theta)\leq r\right)&=\sum_{x\in \mathbb{Z}^4}p_{r-\theta}(O, x)\hat{\mathbb{P}}\left(X_t=O\text{~for some~}t\leq \theta\Big|X_0=x\right)\\
&\leq \sum_{x\in \mathbb{Z}^4}p_{r-\theta}(O, x)\hat{\mathbb{P}}\left(X_t=O\text{~for some~}t<+\infty\Big|X_0=x\right)\\
&=\sum_{x\in \mathbb{Z}^4}p_{r-\theta}(O, x)\frac{\int_0^{+\infty}p_v(x,O)dv}{\int_0^{+\infty}p_v(O, O)dv}\\
&=8\gamma_4\int_0^{+\infty}p_{r-\theta+v}(O, O)dv.
\end{align*}
Hence,
\begin{align*}
\int_{sN}^{tN}\left(\int_{sN}^r\hat{\mathbb{P}}\left(\hat{\tau}_{34}(r, \theta)\leq r\right)d\theta\right)dr
&\leq 8\gamma_4 \int_{sN}^{tN}\left(\int_{sN}^r\left(\int_0^{+\infty}p_{r-\theta+v}(O, O)dv\right)d\theta\right)dr\\
&=8\gamma_4 \int_{sN}^{tN}\left(\int_{0}^{r-sN}\left(\int_0^{+\infty}p_{\theta+v}(O, O)dv\right)d\theta\right)dr\\
&\leq 8\gamma_4 \int_{sN}^{tN}\left(\int_{0}^{(t-s)N}\left(\int_0^{+\infty}p_{\theta+v}(O, O)dv\right)d\theta\right)dr\\
&=8\gamma_4(t-s)N\int_{0}^{(t-s)N}\left(\int_0^{+\infty}p_{\theta+v}(O, O)dv\right)d\theta.
\end{align*}
According to the local central limit theorem of the simple random walk, there exists $K_7<+\infty$ independent of $t$ and $K_8<+\infty$ such that
\[
p_t(O, O)\leq \frac{K_7}{t^2}
\]
for all $t\geq 0$ and $\int_0^{+\infty}p_v(O,O)dv<K_8$. Hence, let $c=\min\{1, (t-s)N\}$, we have
\begin{align*}
&\int_{0}^{(t-s)N}\left(\int_0^{+\infty}p_{\theta+v}(O, O)dv\right)d\theta \\
&=\int_{c}^{(t-s)N}\left(\int_0^{+\infty}p_{\theta+v}(O, O)dv\right)d\theta+\int_{0}^{c}\left(\int_0^{+\infty}p_{\theta+v}(O, O)dv\right)d\theta\\
&\leq K_7\int_c^{(t-s)N}\frac{1}{\theta}d\theta+cK_8\\
&\leq K_7\left(\log(t-s)+\log N-\log c\right) +K_8\leq \log N\left(2K_7+K_8\right)
\end{align*}
for $0\leq s<t\leq T$ and $N\geq \max\{4, T\}$. As a result, \eqref{equ 6.9} holds with $K_6=8\gamma_4\left(2K_7+K_8\right)$ and the proof of Theorem \ref{theorem 2.1 main sample path CLT} in the $d=4$ case is complete.

\qed

\section{The proof of Theorem \ref{theorem 2.1 main sample path CLT}: $d\geq 5$ case}\label{section seven d>=5}
In this section, we prove Theorem \ref{theorem 2.1 main sample path CLT} in the case $d\geq 5$. Hence, throughout this section, we assume that the dimension $d\geq 5$. Our proof follows the same strategy as that given in Section \ref{section six d=4} and hence we only give an outline of this proof to avoid repeating many similar details. We first introduce some notations for later use. For any $d\geq 5, N\geq 1, x\in \mathbb{Z}^d$ and $t\geq 0$, we let $\phi_N(x)$, $\mathcal{H}_N$, $\widetilde{\mathcal{M}}_t^N$ and $\langle \widetilde{\mathcal{M}}^N\rangle_t$ be defined as in Section \ref{section six d=4} except that $4$ in each superscript, which indicates the dimension, is replaced by $d$. Then, we have the following analogue of \eqref{equ d=4 2}. For all $t\geq 0$,
\begin{equation}\label{equ d>=5 1}
\frac{\int_0^{tN}\left(\eta_s(O)-p\right)ds}{\sqrt{N}}
=\frac{\widetilde{\mathcal{M}}_{tN}^N}{\sqrt{N}}-\frac{\mathcal{H}_N(\eta_{tN})}{\sqrt{N}}+\frac{\mathcal{H}_N(\eta_{0})}{\sqrt{N}}
+\frac{\int_0^{tN}\mathcal{H}_N(\eta_s)ds}{\sqrt{N}N}.
\end{equation}
The following lemma is an analogue of Lemma \ref{lemma d=4 1 small error}.
\begin{lemma}\label{lemma d>=5 1 small error}
Let $\eta_0$ be distributed with $\nu_p$. For any given $t\geq 0$,
\[
\lim_{N\rightarrow+\infty}\left(-\frac{\mathcal{H}_N(\eta_{tN})}{\sqrt{N}}+\frac{\mathcal{H}_N(\eta_{0})}{\sqrt{N}}+\frac{\int_0^{tN}\mathcal{H}_N(\eta_s)ds}{\sqrt{N}N}\right)=0
\]
in $L^2$.
\end{lemma}

\proof

As in the proof of Lemma \ref{lemma d=4 1 small error}, we only need to show that
\begin{equation}\label{equ d>=5 3}
\lim_{N\rightarrow+\infty}\frac{1}{N}\left(\sup_{s\geq 0}\mathbb{E}_{\nu_p}\left(\left(\mathcal{H}_N(\eta_s)\right)^2\right)\right)=0.
\end{equation}
According to an analysis similar with that given in the proof of Lemma \ref{lemma d=4 1 small error}, we have
\[
\mathbb{E}_{\nu_p}\left(\left(\mathcal{H}_N(\eta_s)\right)^2\right)
\leq 2dp(1-p)\gamma_d\int_0^{+\infty}N^2\theta e^{-\theta}\left(\int_0^{+\infty}p_{r+N\theta}(O, O)dr\right)d\theta.
\]
When $d\geq 5$, according to the local central limit theorem of the simple random walk, $\int_0^{+\infty}p_r(O, O)dr<+\infty$ and $p_t(O, O)=O(t^{-\frac{d}{2}})$, hence
\begin{align*}
&\int_0^{+\infty}N^2\theta e^{-\theta}\left(\int_0^{+\infty}p_{r+N\theta}(O, O)dr\right)d\theta\\
&=\int_0^{\frac{1}{N}}N^2\theta e^{-\theta}\left(\int_0^{+\infty}p_{r+N\theta}(O, O)dr\right)d\theta
+\int_{\frac{1}{N}}^{+\infty}N^2\theta e^{-\theta}\left(\int_0^{+\infty}p_{r+N\theta}(O, O)dr\right)d\theta\\
&\leq O(1)\int_0^{\frac{1}{N}}N^2\theta e^{-\theta}d\theta+O(1)\int_{\frac{1}{N}}^{+\infty}N^2\theta e^{-\theta}(N\theta)^{1-\frac{d}{2}}d\theta\\
&=N^2O(N^{-2})+O(N^{3-\frac{d}{2}})\int_0^N \theta^{\frac{d}{2}-4}e^{-\frac{1}{\theta}}d\theta\\
&=O(1)+O(N^{3-\frac{d}{2}})\int_0^N \theta^{\frac{d}{2}-4}e^{-\frac{1}{\theta}}d\theta.
\end{align*}
It is easy to check that
\[
O(N^{3-\frac{d}{2}})\int_0^N \theta^{\frac{d}{2}-4}e^{-\frac{1}{\theta}}d\theta=
\begin{cases}
O(\sqrt{N}) & \text{~if~}d=5,\\
O(\log N) & \text{~if~}d=6,\\
O(1) & \text{~if~}d\geq 7.
\end{cases}
\]
In conclusion, $\int_0^{+\infty}N^2\theta e^{-\theta}\left(\int_0^{+\infty}p_{r+N\theta}(O, O)dr\right)d\theta=O(\sqrt{N})$ and hence
\begin{equation}\label{equ 7.3}
\frac{1}{N}\left(\sup_{s\geq 0}\mathbb{E}_{\nu_p}\left(\left(\mathcal{H}_N(\eta_s)\right)^2\right)\right)=O(N^{-\frac{1}{2}}).
\end{equation}
Equation \eqref{equ d>=5 3} follows from \eqref{equ 7.3} and the proof is complete.

\qed

The following lemma is an analogue of Lemma \ref{lemma d=4 finite dimension}, which shows that any finite dimensional distribution of $\left\{\frac{1}{\sqrt{N}}\int_0^{tN}\left(\eta_s(O)-p\right)ds\right\}_{0\leq t\leq T}$ converges weakly to the corresponding finite dimensional distribution of $\{C_dW_t\}_{0\leq t\leq T}$.

\begin{lemma}\label{lemma d>=5 finite dimension}
Let $\eta_0$ be distributed with $\nu_p$, then for each integer $m\geq 1$ and $0\leq t_1<t_2<\ldots<t_m\leq T$,
\[
\frac{1}{\sqrt{N}}\left(\int_0^{t_1N}\left(\eta_s(O)-p\right)ds, \int_0^{t_2N}\left(\eta_s(O)-p\right)ds,\ldots, \int_0^{t_mN}\left(\eta_s(O)-p\right)ds\right)
\]
converges weakly, as $N\rightarrow+\infty$, to $C_d\left(W_{t_1}, W_{t_2}, \ldots, W_{t_m}\right)$.
\end{lemma}

As in Section \ref{section six d=4}, we require the following analogue of Lemma \ref{lemma d=4 convergence of mathcalH} as a preliminary of the proof of Lemma \ref{lemma d>=5 finite dimension}.

\begin{lemma}\label{lemma d>=5 convergence of mathcalH} For $d\geq 5$, we have
\[
\lim_{N\rightarrow+\infty}\sum_{x\in \mathbb{Z}^d}\phi_N^2(x)=\int_0^{+\infty}sp_s(O, O)ds.
\]
\end{lemma}

\proof

As we have recalled in Section \ref{section two}, when $d\geq 5$, $\int_0^{+\infty}sp_s(O, O)ds<+\infty$. Hence,
\begin{align*}
\lim_{N\rightarrow+\infty}\sum_{x\in \mathbb{Z}^d}\phi_N^2(x)&=\lim_{N\rightarrow+\infty}\int_0^{+\infty}se^{-\frac{s}{N}}p_s(O, O)ds\\
&=\int_0^{+\infty}sp_s(O, O)ds
\end{align*}
and the proof is complete.

\qed

Now, we prove Lemma \ref{lemma d>=5 finite dimension}.

\proof[Proof of Lemma \ref{lemma d>=5 finite dimension}]

By \eqref{equ d>=5 1} and Lemma \ref{lemma d>=5 1 small error}, we only need to check that
\begin{equation}\label{equ d>=5 convergence of quadratic}
\lim_{N\rightarrow+\infty}\frac{1}{N}\langle \widetilde{\mathcal{M}}^N\rangle_{tN}=C_d^2 t
\end{equation}
in $L^2$ for any $t\geq 0$. As in the proof of Lemma \ref{lemma d>=5 finite dimension}, we have
\[
\frac{1}{N}\langle \widetilde{\mathcal{M}}^N\rangle_{tN}=\sum_{x\in \mathbb{Z}^d}\sum_{y\sim x}\int_0^{t}\phi_N^2(x)\left(\eta_{sN}(y)-\eta_{sN}(x)\right)^2ds
\]
by utilizing the Dynkin's martingale formula. Then, Lemma \ref{lemma d>=5 convergence of mathcalH} and \eqref{equ duality minus square} ensures that
\[
\lim_{N\rightarrow+\infty}\mathbb{E}_{\nu_p}\left(\frac{1}{N}\langle \widetilde{\mathcal{M}}^N\rangle_{tN}\right)=C_d^2t.
\]
Hence, to complete the check of \eqref{equ d>=5 convergence of quadratic}, we only need to show that
\begin{equation}\label{equ 7.5}
\lim_{N\rightarrow+\infty}\frac{1}{N^2}{\rm Var}_{\nu_p}\left(\langle \widetilde{\mathcal{M}}^N\rangle_{tN}\right)=0.
\end{equation}
As in the proofs of Lemmas \ref{lemma 4.3} and \ref{lemma d=4 finite dimension}, Equation \eqref{equ 7.5} follows from a $d\geq 5$ version of \eqref{equ 4.7}. As we have explained in the proof of Lemma \ref{lemma d=4 finite dimension}, \eqref{equ 4.7} can be extended to all $d\geq 3$ cases since simple random walks on $\mathbb{Z}^d$ with $d\geq 3$ are transient. Since \eqref{equ 7.5} holds, the proof is complete.

\qed

At last, we prove Theorem \ref{theorem 2.1 main sample path CLT} in the case $d\geq 5$.

\proof[Proof of Theorem \ref{theorem 2.1 main sample path CLT} in case $d\geq 5$]

According to an analysis similar with that given in the proof of the case $d=4$, we only need to show that there exists $K_{9}<+\infty$ independent of $N, s, t$ such that
\begin{equation}\label{equ 7.6}
\frac{1}{N}\int_{sN}^{tN}\left(\int_{sN}^r\hat{\mathbb{P}}\left(\hat{\tau}_{34}(r, \theta)\leq r\right)d\theta\right)dr\leq K_9(t-s)
\end{equation}
for all $N\geq 1$ and $0\leq s<t\leq T$, where $\hat{\tau}_{34}(r, \theta)$ is the $\hat{\tau}_{34}$ defined in Section \ref{section three} with $t_3=\theta, t_4=r$ and $(x_3, x_4)=(O, O)$. According to a calculation similar with that given in the proof of the case $d=4$, we have
\begin{align*}
\int_{sN}^{tN}\left(\int_{sN}^r\hat{\mathbb{P}}\left(\hat{\tau}_{34}(r, \theta)\leq r\right)d\theta\right)dr&
\leq 2d\gamma_d(t-s)N\int_{0}^{(t-s)N}\left(\int_0^{+\infty}p_{\theta+v}(O, O)dv\right)d\theta\\
&\leq 2d\gamma_d(t-s)N\int_{0}^{+\infty}\left(\int_0^{+\infty}p_{\theta+v}(O, O)dv\right)d\theta\\
&= 2d\gamma_d(t-s)N\int_0^{+\infty}\theta p_\theta(O, O)d\theta.
\end{align*}
Since $\int_0^{+\infty}\theta p_\theta(O, O)d\theta<+\infty$ when $d\geq 5$, \eqref{equ 7.6} holds with $K_9=2d\gamma_d\int_0^{+\infty}\theta p_\theta(O, O)d\theta$ and the proof is complete.

\qed

\section{Comments on Conjecture \ref{conjecture in d=2 case}}\label{section d=2 conjecture comment}

In this section, we explain why we believe Conjecture \ref{conjecture in d=2 case} holds and what we lack to prove this conjecture. Let $N_t^{x,y}$, $\hat{N}_t^{x,y}$, $v(t,x)$, $V_s^t$, $\mathcal{M}_s^t$ be defined as in Section \ref{section four d=3} except that $3$ in each superscript indicating the dimension is replaced by $2$, then \eqref{equ 4.2} has a $d=2$ version and
\begin{equation}\label{equ 5.1}
\frac{1}{\frac{N}{\sqrt{\log N}}}\int_0^{tN}(\eta_s(O)-p)ds=\frac{1}{\frac{N}{\sqrt{\log N}}}\mathcal{M}_{tN}^{tN}+\frac{1}{\frac{N}{\sqrt{\log N}}}V_0^{tN}(\eta_0).
\end{equation}

We have the following $d=2$ version of Lemma \ref{lemma small error}.
\begin{lemma}\label{lemma small error d=2}
For any given $t>0$, if $\eta_0$ is distributed with $\nu_p$, then
\[
\lim_{N\rightarrow+\infty}\frac{1}{\frac{N}{\sqrt{\log N}}}V_0^{tN}(\eta_0)=0
\]
in $L^2$.
\end{lemma}

By utilizing the fact that $p_t(O, O)=O(t^{-1})$ when $d=2$, Lemma \ref{lemma small error d=2} follows from an analysis similar with that leading to Lemma \ref{lemma small error}.

For any $u\in \mathbb{R}^2$ and integer $N\geq 1$, let $u_N$ be defined as in Section \ref{section four d=3} except that $3$ in each superscript indicating the dimension is replaced by $2$. For each $N\geq 1$, we denote by $\hat{\mathcal{Y}}^N$ the random measure on $[0, T]\times \mathbb{R}^2$ such that
\begin{align*}
&\hat{\mathcal{Y}}^N(H)=\sqrt{\log N}\int_{\mathbb{R}^2}\sum_{y:y\sim \sqrt{N}u_N}\left(\int_0^TH(s,u)\left(\eta_{Ns-}(y)-\eta_{Ns-}(\sqrt{N}u_N)\right)d\hat{N}_{Ns}^{\sqrt{N}u_N,y}\right)du\\
&=\sqrt{\log N}\sum_{x\in \mathbb{Z}^2}\sum_{y:y\sim x}\int_0^T\left(\int_{\frac{x}{\sqrt{N}}+(-\frac{1}{2\sqrt{N}}, \frac{1}{2\sqrt{N}}]^2} H(s,u)du\right)\left(\eta_{Ns-}(y)-\eta_{Ns-}(x)\right)d\hat{N}_{Ns}^{x,y}
\end{align*}
for any $H\in C_c\left([0, T]\times \mathbb{R}^2\right)$. We denote by $\hat{\mathcal{Y}}$ the Gaussian time-space white noise on $[0, T]\times \mathbb{R}^2$ such that
$\hat{\mathcal{Y}}(H)$ follows the normal distribution with mean zero and variance
\[
8p(1-p)\pi\int_0^T\int_{\mathbb{R}^2}H^2(s,u)dsdu
\]
for any  $H\in C_c\left([0, T]\times \mathbb{R}^2\right)$.

For each $N\geq 1$ and $0\leq t\leq T$, we define
\[
\hat{b}_t^N(s,u)=\sum_{x\in \mathbb{Z}^2}v\left(N(t-s), x\right)1_{\{u\in \frac{x}{\sqrt{N}}+(-\frac{1}{2\sqrt{N}}, \frac{1}{2\sqrt{N}}]^2, s\leq t\}}
\]
for any $(s, u)\in [0, T]\times \mathbb{R}^2$. Then, $\hat{b}_t^N$ converges, as $N\rightarrow +\infty$, to $\hat{b}_t$ uniformly and in $L^2$ on $[0, T]\times \mathbb{R}^2$, where
\[
\hat{b}_t(s,u)=\int_0^{t-s}\frac{1}{4\pi r}e^{-\frac{\|u\|_2^2}{4r}}dr1_{\{s\leq t\}}.
\]
By the $d=2$ version of \eqref{equ 4.2}, we have $\frac{1}{\frac{N}{\sqrt{\log N}}}\mathcal{M}_{tN}^{tN}=\hat{\mathcal{Y}}^N(\hat{b}_t^N)$. 

According to the fact that $\lim_{t\rightarrow+\infty}(\log t)\hat{\mathbb{P}}\left(\tau_{xy}>t\right)=\pi$ for $x,y\in \mathbb{Z}^2$ such that $x\sim y$ (see (1.5) of \cite{Cox1983}) and $\lim_{t\rightarrow+\infty}tp_{t,2}(O, O)=\frac{1}{4\pi}$ (see Chapter 2 of \cite{Lawler2010}), it is not difficult to check the following analogue of \eqref{equ 4.11}.
For any $0\leq s<t\leq T$,
\begin{align}\label{equ d=2 4.11 analogue}
&\lim_{N\rightarrow+\infty}{\rm Cov}\left(\frac{\sqrt{\log N}}{N}\mathcal{M}_{sN}^{sN}, \frac{\sqrt{\log N}}{N}\mathcal{M}_{tN}^{tN}\right)\\
&=C_2^2\left(\frac{(t+s)^2}{4}\log(t+s)+\frac{(t-s)^2}{4}\log(t-s)-\frac{s^2\log s}{2}-\frac{t^2\log t}{2}\right). \notag
\end{align}

According to \eqref{equ d=2 4.11 analogue} and an analysis similar with that given in the proof of Theorem \ref{theorem 2.1 main sample path CLT} in the case $d=3$, it is not difficult to check that Conjecture \ref{conjecture in d=2 case} is a corollary of the following conjecture, which is a $d=2$ version of Lemma \ref{lemma 4.4}. 

\begin{conjecture}\label{lemma 5.2}
Let $\eta_0$ be distributed with $\nu_p$, then for any $H\in C_c\left([0, T]\times \mathbb{R}^2\right)$, $\hat{\mathcal{Y}}^N(H)$ converges weakly to $\hat{\mathcal{Y}}(H)$ as $N\rightarrow+\infty$ and furthermore,
\[
\lim_{N\rightarrow+\infty}\mathbb{E}_{\nu_p}\left(\left(\hat{\mathcal{Y}}^N(H)\right)^2\right)=\mathbb{E}\left(\left(\hat{\mathcal{Y}}(H)\right)^2\right).
\]
\end{conjecture}

It is natural to attempt to prove Conjecture \ref{lemma 5.2} through an analysis similar with that given in the proof of Lemma \ref{lemma 4.4}. According to the fact that $\lim_{t\rightarrow+\infty}(\log t)\hat{\mathbb{P}}\left(\tau_{xy}>t\right)=\pi$, we can check an analogue of \eqref{equ 4.5 two}. However, currently we can not check a $d=2$ version of \eqref{equ 4.6}, since our approach given in the $d=3$ case relies heavily on the transience of the simple random walk. We will work on the $d=2$ case as a further investigation. 

\quad

\textbf{Acknowledgments.}
The author is grateful to financial
supports from the Fundamental Research Funds for the Central Universities with grant number 2022JBMC039 and the National Natural Science Foundation of China with grant number 12371142.

{}
\end{document}